\documentclass[12pt]{amsart}
\usepackage{Preamble}

{
\newtheorem*{t1}{Theorem \ref{module}}
\newtheorem*{t2}{Theorem \ref{Thm}}
}

\title{A Quantum $H^*(T)$-module via Quasimap Invariants}
\author{Jae Hwang Lee}

\begin{document}

\begin{abstract}
    For $X$ a smooth projective variety, the quantum cohomology ring $QH^*(X)$ is a deformation of the usual cohomology ring $H^*(X)$, where the product structure is modified to incorporate quantum corrections. These correction terms are defined using Gromov--Witten invariants. When $X$ is toric with the geometric quotient description $V\git T$, the cohomology ring $H^*(V\git T)$ also has the structure of a quantum $H^*(T)$-module. In this paper, we give a new deformation using quasimap invariants with a light point. This defines $H^*(T)$-module structure on $H^*(X)$ through a modified version of the WDVV equations. Using the Atiyah--Bott localization theorem, we explicitly compute this structure for the Hirzebruch surface of type 2. We conjecture that this new quantum module structure is isomorphic to the natural module structure of the Batyrev ring for a semipositive toric variety.
\end{abstract}
\maketitle
\tableofcontents
\setcounter{tocdepth}{1} 

\section{Introduction}
Quantum cohomology is a central object within the interests of both mathematics and physics, as it is related to string theory and mirror symmetry \cite{morrison}. The idea of quantum cohomology first appeared in physics \cite{lerche}. The first mathematical construction was given in terms of symplectic geometry for semi-positive symplectic manifolds \cite{mcduff,ruan}.

Quantum cohomology has a ring structure whose product is called the \textit{quantum product}. The quantum product is a deformation of the product of the ordinary cohomology. In algebraic geometry side, it is defined by using Gromov--Witten invariants \cite{coxmirror,mirrsymm,kont} via the moduli space of stable maps. One notable application of the quantum product structure, especially \textit{associativity}, has been shown in \cite{kont} by verifying the prediction given in \cite{candelas} on the number of rational curves of degree 4 on a quintic 3-fold.

However, it is not easy to compute quantum cohomology in general. When the space $X$ is a smooth Fano toric variety, the quantum cohomology ring $QH^*(X)$ agrees with the Batyrev ring $\text{Bat}^*(X)$, defined in \cite{batyrev} (see \cite[Example 8.1.2.2 or Example 11.2.5.2]{coxmirror} or an extension of \cite{ruanqin}). On the other hand, the Hirzebruch surface of type 2, say $\F_2$, which is not Fano, but \textit{semipositive}, shows a failure of such an equality \cite[Example 11.2.5.2]{coxmirror}.
\vspace{0.3cm}

In this paper, we define a new operation with different quantum deformations (Definition \ref{21oper}) of the product of the ordinary cohomology using $2|1$-quasimap invariants, instead of $3$-pointed GW invariants. The moduli space of $m|k$-pointed quasimaps were introduced in \cite{bigI} (see Definition \ref{defquasimap}). Here, $m|k$ means the number of heavy and light markings, respectively, where heavy markings are the ordinary ones and light markings are infinitesimally weighted ones \cite{lho}. This operation defines not a product structure, but a module structure; the operation satisfies the compatibility \eqref{wdvvintro} with the cup product in a cohomology ring. This is exactly an analogue of associativity of the quantum product, or Witten-Dijkgraaf-Verlinde-Verlinde (WDVV) equations.

Let $X_\Sigma$ be a smooth projective toric variety whose toric geometric quotient description (Theorem \ref{toricquo}) is given by $V \git T$, where $V$ is a finite dimensional $\C$-vector space and $T$ is a complex torus. Denote $H^*(T)$ the group cohomology of $T$. Then, the following is the main theorem which defines the quantum $H^*(T)$-module structure on $H^*(X_\Sigma)$:
\begin{t1}[Quantum $H^*(T)$-module structure]
For $\xi, \zeta \in H^*(T)$ and $\phi \in H^*(X_\Sigma)$,
\begin{equation}
\label{wdvvintro}
\xi \star (\zeta \star \phi) = (\xi  \cdot \zeta ) \star \phi .
\end{equation}
\end{t1}
We call this structure the \textit{(small) quantum $H^*(T)$-module structure}.

An explicit computation for the quantum module structure of $\F_2$ will be given in Section \ref{sec_compute} by applying the Atiyah--Bott localization theorem. Suppose that $\Pic(\F_2)\simeq \Z D_2 \oplus D_4$, where $D_2$ and $D_4$ are the torus-invariant divisors such that $D_2\cdot D_2=0$ and $D_4\cdot D_2=-2pt$. The geometric quotient construction of $\F_2$ is $\C^4 \git (\C^*)^2$. The group cohomology of $(\C^*)^2$ is given by $\C[q_2,q_4]$. A full description of the quantum module structure of $\F_2$ is given as follows:
\begin{t2}
The quantum $H^*((\C^*)^2)$-module structure for $\F_2$ is given by the following
\begin{align*}
    \sigma_2 \star 1 &= D_2-\frac 12 f(q_4)D_4 &\quad \sigma_4 \star 1&= (1+f(q_4))D_4 \\
    \sigma_2 \star D_2 &=q_2q_4(1+f(q_4)) -\frac{1}{2}f(q_4)pt &\quad \sigma_4 \star D_2&=-\frac{1}{2}q_2f(q_4) + (1+f(q_4))pt\\
    \sigma_2 \star D_4 &=-2q_2q_4(1+f(q_4)) + (1+f(q_4))pt &\quad \sigma_4 \star D_4&=q_2(1 + f(q_4))-2(1+f(q_4))pt\\
    \sigma_2 \star pt &= q_2q_4(1+f(q_4))D_4 &\quad \sigma_4 \star pt&= q_2 D_2 - \frac{1}{2}q_2(1+f(q_4))D_4,
\end{align*}
where $f(z)=\sum_{d\geq 1}\binom{2d}{d}z^d = \frac{1}{\sqrt{1-4z}}-1$.
\end{t2}

An interesting observation of this result is that the quantum module structure of $\F_2$ coincides with the Batyrev ring of $\F_2$ regarded as a module. In other words, we found a geometric interpretation of the Batyrev ring of $\F_2$ through $2|1$-quasimap invariants. Based on this evidence, we conjecture the following:
\begin{Conj}
\label{conjbatyrev}
    \textit{For a smooth semipositive toric variety $V \git T$, the quantum $H^*(T)$-module structure of $V \git T$ coincides with a natural module structure of the Batyrev ring of $V \git T$.}
\end{Conj}

To prove Theorem \ref{module}, which is an analogue of WDVV equations, a type of a splitting lemma as in \cite{abramo,lucasplit,kock} is required. We need to deal with the virtual fundamental classes of the quasimap moduli spaces and the diagonal pullback. Instead of directly using the perfect obstruction theory as in \cite{lucasplit}, we use the localized top Chern classes \cite[ch14.1]{fulton}, which is more elementary notion. There is a construction of the moduli space of quasimaps with light points as a zero locus of a section of a vector bundle on a smooth DM stack \cite{toricquasi,gitquasi}. This global model gives rise to the virtual fundamental classes as the localized top Chern, which agrees with the Behrend--Fantachi version in \cite{behfan} (see \cite{mark}). 

In our localization computation for the quantum module structure of $\F_2$, there are some technical key features that we would like to highlight:
\begin{enumerate}
    \item The module structure in Theorem \ref{module} allows us to assume that the degree of the insertion from the light marking is one.
    \item In general, the map forgetting a heavy point does not define a universal curve in quasimap case. However, the map forgetting a light point gives rise to a universal curve \cite{lho}. This allow us to have the divisor equation as in \cite{jinzenji}. Thus, the computation of the module structure boils down to the computation of all possible 2-pointed quasimap invariants.
    \item Having 2 heavy markings gives rise to a chain of $\PP^1$'s for the source curve of a quasimap.
    \item When applying the localization theorem, a contribution of one fixed locus can be expressed as a fraction whose numerator is a homogeneous polynomial in $V$ and $W$, where $V$ and $W$ are weights given in \eqref{vandw}, and the denominator is $W^N$ for some $N$. This allows us to look at particular type of fixed loci, called the necessary fixed loci (see the paragraph above Corollary \ref{correducedloci}, and Figure \ref{reducedloci}, and Figure \ref{reducedloci2}).
\end{enumerate}
\vspace{0.3cm}

This paper is organized as follows. In Section \ref{sec_pre}, we give a preliminaries on the quotient construction of a toric variety and recall the definition of the moduli space of stable toric quasimaps with light points. Such a moduli space of quasimaps will be constructed with a global embedding in Section \ref{sec_construction}. The definition of the new quantum deformation using these moduli spaces is given in Section \ref{sec_splitting}. We also prove Splitting Lemma \ref{split}. In Section \ref{sec_compute}, we elaborate our computation for the quantum module structure of the Hirzebruch surface of type 2, using the Atiyah--Bott localization theorem. In Section \ref{sec_conj}, we verify that the quantum module structure of $\F_2$ agrees with the Batyrev ring realized as a module. 

\textbf{Acknowledgement.} I have benefited from conversations with Jeongseok Oh and Woonam Lim on the matter of virtual fundamental classes and localization computations, respectively. Special thanks to Qaasim Shafi for the suggestion of Luca Battistella's dissertation \cite{lucathesis}, from which I was able to connect this research with the Batyrev rings. Last but not least, I am thankful for the many helpful discussions with Renzo Cavalieri and Mark Shoemaker. This project was partially supported by Renzo Cavalieri's NSF DMS 2100962, Mark Shoemaker's NSF Grant 1708104, and fully supported by the Department of Mathematics Summer Research Fellowship for 2023 at Colorado State University. I gratefully acknowledge financial support from Douglas Ortego.

\section{Moduli Spaces of Stable Toric $m|k$-pointed Quasimaps to $X_\Sigma$}
\label{sec_pre}
Through out this paper, our base field is the complex numbers $\C$. Let $M$ be a $\Z$-lattice, and $N$ the dual lattice, and $\Sigma \subseteq N_\R:=N \otimes_\Z \R$ a smooth complete fan. Write $X_\Sigma$ for the corresponding smooth projective toric variety with the torus $N \otimes_\Z \C^*$.

\subsection{Geometric Quotients for Toric Varieties} The toric variety $X_\Sigma$ can be expressed as a geometric quotient using data of the fan $\Sigma$. A \textit{primitive collection} is a subset $P$ of rays in $\Sigma(1)$ such that
\begin{enumerate}[label=\roman*)]
    \item $P$ is not contained in any $\sigma\in \Sigma$;
    \item every proper subset of $P$ is contained in $\sigma(1)$ for some $\sigma \in \Sigma$.
\end{enumerate} 
Define 
\[Z_\Sigma := \bigcup_{P: \text{ a primitive collection}} \V(x_\rho \mid \rho \in P) \subseteq \C^{\Sigma(1)}  ,\]
where $\rho$ is chosen to be the minimal generator of the ray. This is the \textit{irrelevant subset} to the fan $\Sigma$.

Since the fan $\Sigma$ is complete, we have an exact sequence
\begin{equation}
\label{seqclassgroup}
    0 \longrightarrow M \longrightarrow \Z^{\Sigma(1)}\longrightarrow \text{Cl}(X_\Sigma) \longrightarrow 0,
\end{equation}
where $m\in M$ goes to $\sum_{\rho \in \Sigma(1)} \langle m, \rho  \rangle$. Smoothness of $X_\Sigma$ allows us to identify the class group $\text{Cl}(X_\Sigma)$ with the Picard group $\text{Pic}(X_\Sigma)$. Denote the matrix of the map $\Z^{\Sigma(1)} \rightarrow \text{Pic}(X_\Sigma)$ by $(a_{i\rho})$ where $\rho=1,2,\ldots, n$ and $i=1,2,\ldots, r$ with $n:=|\Sigma(1)|$ and $r:= \rank \Pic(X_\Sigma)$. Applying $\Hom_\Z(-,\C^*)$, we obtain the exact sequence
\[0 \longrightarrow (\C^*)^r \longrightarrow (\C^*)^{\Sigma(1)} \longrightarrow  N \otimes_\Z \C^* \longrightarrow 0.\]
Thus, the torus $(\C^*)^r$ is acting on $\C^{\Sigma(1)}$ through the componentwise multiplication, which can be represented by the transposition of the $r\times n$ matrix  $(a_{i\rho})$. In this case, the geometric quotient associated to $X_\Sigma$ is given as follows from \cite[Theorem 5.1.11]{coxtoric}.
\begin{Thm}
\label{toricquo}
Given a smooth complete fan $\Sigma$, there is a natural isomorphism between the corresponding toric variety and the geometric quotient \[X_\Sigma \simeq \C^{\Sigma(1)} \backslash Z_\Sigma ~\git~ (\C^*)^r.\]
\end{Thm}

\subsection{The moduli space of stable toric $m|k$-quasimaps} We recall the definition of the moduli space of stable toric quasimaps with $m$ heavy points and $k$ light points to a smooth projective toric variety $X_\Sigma$. For details, we refer to \cite{toricquasi} for a construction of the moduli space of stable toric quasimaps with heavy markings and \cite{bigI} for the one with light markings.

Choose $\alpha_\rho\in \Z$ so that for a line bundle $\calO_{X_\Sigma}(1):=\otimes_{\rho \in \Sigma(1)}\calO(D_\rho)^{\otimes \alpha_\rho}$ on $X_\Sigma$ to be ample, where $D_\rho$ is the torus invariant divisor on $X_\Sigma$ corresponding to $\rho$. 
\begin{Def}
\label{defquasimap}
A \textit{stable toric $m|k$-pointed quasimap} to $X_\Sigma$ of genus $g$ is the data \[((C;x_1,\ldots,x_m;y_1,\ldots,y_k),\{L_\rho \}_{\rho\in \Sigma(1)}, \{s_\rho\}_{\rho \in \Sigma(1)} , \{\phi_m\}_{m\in M} )\]
where
\begin{itemize}
    \item $C$ is a connected, at most nodal, projective curve of genus $g$,
    \item $\{x_1,\ldots, x_m,y_1,\ldots,y_k\}$ are nonsingular marked points,
    \item $\{x_1,\ldots, x_m\}$ are distinct and disjoint from $\{y_1,\ldots,y_k\}$,
    \item $L_\rho$ are line bundles on $C$,
    \item $s_\rho \in \Gamma (C, L_\rho)$ are global sections,
    \item (compatibility) the trivializations $\phi_m: \otimes_\rho L_\rho^{\langle m,\rho \rangle } \rightarrow \mathcal{O}_{C}$ are isomorphisms satisfying $\phi_m \otimes \phi_{m'}=\phi_{m+m'} $ for all $m,m' \in M$,
\end{itemize} satisfying
\begin{enumerate}
    \item (nondegeneracy) there is a finite (possibly empty) set of smooth points $B \subset C$, disjoint from $\{x_1,\ldots, x_m\}\subset C$, such that for every $z \in C \backslash B$, there exists a maximal cone $\sigma \in \Sigma_{max}$ such that $u_\rho(z)  \neq 0,~ \forall \rho \not\subset \sigma$,
    \item (stability) $\omega_C (x_1+\cdots+x_m+\epsilon(y_1+\cdots+y_k) ) \otimes \mathcal{L}^{\epsilon}$ is ample for every rational number $\epsilon > 0$, where $\mathcal{L} := \otimes_{\rho \in \Sigma(1)} L_\rho^{\otimes \alpha_\rho}$.
\end{enumerate}
\end{Def}

Denote a stable quasimap by $(C;\underline{x};\underline{y},\underline{L},\underline{s},\underline{\phi})$. Two stable quasimaps $(C;\underline{x};\underline{y},\underline{L},\underline{s},\underline{\phi})$ and $(C';\underline{x'};\underline{y'},\underline{L'},\underline{s'},\underline{\phi'})$ are \textit{isomorphic} if there exists \[(f:C\rightarrow C', \{\theta_\rho : L_\rho \rightarrow f^*(L'_\rho )  \}_{\rho\in \Sigma(1)} ), \] where $f$ and $\theta_\rho$ are isomorphism such that 
\[f(x_i)=x'_i,~f(y_i)=y'_i,~\theta_\rho (s_\rho)=f^*(s'_\rho),~\phi_m = f^*(\phi'_m) \circ (\otimes_\rho \theta_{\rho}^{\langle m,\rho  \rangle}  ).  \]

\begin{Def}
Given a stable toric quasimap $(C;\underline{x};\underline{y},\underline{L},\underline{s},\underline{\phi})$,
the map $\Pic(X_\Sigma) \rightarrow \Pic(C)$ sending $\mathcal{O}_{X_\Sigma}(D_\rho) \mapsto L_\rho$ is a well-defined homomorphism because of the compatibility condition of the trivializations $\phi_m$. Composing with the degree map to $\Z$, there is a $\Z$-linear homomorphism from $\Pic(X_\Sigma)$ to $\Z$. By Poincar\'{e} duality, a perfect pairing between $H_2(X_\Sigma,\Z)$ and $H^2(X_\Sigma,\Z) \simeq \Pic(X_\Sigma)$, gives rise to a unique class $\beta \in H_2(X_\Sigma, \Z)$ determined by
\[ \beta \cdot D_\rho = \deg L_\rho, \]
for all $\rho \in \Sigma(1)$.
This $\beta$ is called the \textbf{degree} of the stable quasimap. Denote the degree restricted to a component $C'$ of the source curve of a quasimap by $\beta_{C'}$.
\end{Def}

\begin{Rem}
\begin{enumerate}
    \item We call $x_i$ a \textit{heavy marking} and $y_j$ a \textit{light marking}.
    \item From the degree of the log-canonical bundle $\omega_C(\sum x_i + \epsilon\sum y_j)\otimes \mathcal{L}^\epsilon $, we obtain $2g-2+m\geq 0$.
    \item On a rational component $C'$ of $C$, if $\beta_{C'}> 0$, then $C'$ must have at least two special points, i.e., a heavy marking or a node. When $\beta_{C'}=0$, there are at least three special points, or at most one of them can be replaced by a light point.
    \item On a genus one component, there is at least one special point or a light marking, otherwise the line bundle $\mathcal{L}$ restricted to the component must be of positive degree.
    \item The subset $B \subset C$ is the set of \textit{base points} of a quasimap. Away from each base point, the sections of a quasimap defines a map to $X_\Sigma$.
    \item Observe that the nondegeneracy condition is only related to heavy markings. Thus, light markings can collide with base points.
    \item We do not have any \textit{rational tails}, i.e., a component without any markings, since a quasimap with such a component is not stable. It follows that the number of components is finite, because the number of markings is finite.  
\end{enumerate}
\end{Rem}
\begin{Rem}
    \label{remreducedquasi}
    The compatibility condition of trivializations $\phi_m$ in Definition \ref{defquasimap} can be dropped, since it can be recovered from the matrix $(a_{i\rho})$ of the map $\Z^{\Sigma(1)} \rightarrow \text{Pic}(X_\Sigma)$ in \ref{seqclassgroup}. Thus, choosing an integral basis for $\text{Pic}(X_\Sigma)$, say $\{\mathcal{P}_1,\ldots,\mathcal{P}_r\}$, the data of a degree $\beta$ stable toric quasimap $(C;\underline{x};\underline{y},\underline{L},\underline{s},\underline{\phi})$ is equivalent to the data
    \[(C;\underline{x};\underline{y},\underline{P},\underline{s}),\]
    where $P_i$ are line bundles on $C$ with $\deg P_i= \int_\beta c_1(\mathcal{P}_i)$, and $L_\rho=\bigotimes_{i=1}^{r} P_i^{\otimes a_{i\rho}}$. In Section \ref{sec_construction}, we will use this equivalent description for quasimaps.
    
\end{Rem}
The following moduli space was constructed in \cite{toricquasi,bigI}.
\begin{Def} Fix $g, m, k \geq 0$, $\beta \in H_2(X_\Sigma)$. The \textbf{moduli space of degree $\beta$ stable toric quasimaps} to $X_\Sigma$ is the moduli stack parametrizing isomorphism classes of families of stable toric $m|k$-quasimaps of degree $\beta$. Denote it by $Q_{g,m|k}(X_\Sigma,\beta)$.
\end{Def}

\subsection{Stack quotients} 


Due to the presence of base points, a quasimap does not define a map to the toric variety $X_\Sigma\simeq \C^{\Sigma} \backslash Z_\Sigma \git (\C^*)^r$. Instead, the natural target of a quasimap is the stack quotient $[V/T]$, where $V:= \C^\Sigma$ and $T:=(\C^*)^r$. Note that $X_\Sigma$ is an open substack of $[V/T]$.

The cohomology of the stack quotient is given by
\[H^*([V / T]) = H^*([\text{pt} / T]) = H_{T}^{*}(\text{pt}) = H^*(T) =\C[\sigma_1,\ldots, \sigma_r],\] where $H^*(T)$ is the group cohomology of $T$. The variables $\sigma_j$ have a geometric interpretation. Suppose that we have an integral basis $\{\mathcal{P}_{j}\}_{j=1}^{r}$ of $\text{Pic}(X_\Sigma)$, where $\mathcal{P}_{j}$ correspond to $\calO_{X_\Sigma}(D_j)$. The $j$th column of the $r\times n$ action matrix $(a_{i\rho})$ from the toric quotient defines a $T$-equivariant line bundle $\C \rightarrow $pt. Then, the equivariant Euler class of this bundle gives $\sigma_i \in H_{T}^{*}(\text{pt})$.

In \cite{bigI}, the moduli space of quasimaps with light points $Q_{g,m|k}(X_\Sigma,\beta)$ was identified with the following moduli space of quasimaps to a stack quotient without any light points
\[Q_{g,m|0}([\C^n / (\C^*)^r] \times [\C / \C^*]^{k}, (\beta,1,\ldots,1)) .\]
In \cite{gitquasi}, they showed that such a moduli space has a perfect obstruction theory to define the virtual fundamental class.

\subsection{Evaluation maps} For a $i$th heavy marking, the evaluation map \[ev_i: Q_{g,m|k}(X_\Sigma,\beta) \rightarrow X_\Sigma\] is well-defined since heavy markings are distinct from base points and light markings, so that sections of a quasimap defines a map to $X_\Sigma$.

In contrast, since light markings can collide with base points, sections of a quasimap might not define a map to $X_\Sigma$ to define an evaluation map. In this case, the stack quotient $[V / T]$ can be used as the target of an evaluation map at each light marking, since base points can land on the complement of $X_\Sigma$ in $[V / T]$. Thus, the evaluation map at the $j$th light marking is given as follows:
\[\hat{ev}_j: Q_{g,m|k}(X_\Sigma,\beta) \rightarrow [V / T].\]

For a fixed genus $g$ and a degree $\beta$, the $m|k$-pointed quasimap invariants are defined as follows:
\begin{Def}
For $\phi_1,\ldots,\phi_m \in H^*(X_\Sigma)$ and $\xi_1,\ldots,\xi_k \in H^*([V / T])$, an \textit{$m|k$-pointed quasimap invariant} is defined by
\[\langle \phi_1,\ldots,\phi_m \: | \xi_1,\ldots,\xi_k \rangle_{g,m|k,\beta}:= \int_{[Q_{g,m|k}(X_\Sigma,\beta)]^{\text{vir}}} \prod_{i=1}^{m} ev_i^*(\phi_i) \prod_{j=1}^{k} \hat{ev}_j^*(\xi_j). \]
\end{Def}

\section{Global Construction}
\label{sec_construction}
The quasimap moduli spaces $Q_{g,m|k}(X_\Sigma,\beta)$ have a global description. In other words, it can be realized as a stack-theoretic zero locus $Z(s)$ of a section $s$ of a vector bundle $E$ over a smooth Deligne-Mumford stack $B$. In this case, the virtual class of the moduli space is given by the localized Euler class. We briefly review the localized Euler class in \cite[\S14.1]{fulton}, and construct the virtual class by giving the global setting for our moduli space. A detailed construction is given in \cite{toricquasi}.

\subsection{Localized Euler classes}\label{loceuler} For a fiber square of schemes
\[\begin{tikzcd}
X' \arrow[d, "g"'] \arrow[r, "j"] & Y' \arrow[d, "f"] \\
X \arrow[r, "i"]                    & Y               
\end{tikzcd}\]
with $i$ a regular embedding of codimension $d$, there exists an induced homomorphism between the Chows groups
\[i^!: A_* Y' \rightarrow A_{*} X',\] 
which is called the \textbf{Gysin homomorphism} \cite[\S6.2]{fulton}.

Let $E$ be a rank $e$ vector bundle over a pure $n$-dimensional scheme $X$, and $s$ is a section of $E$. Denote the zero scheme by $Z(s)$. Consider the fiber square
\[\begin{tikzcd}
Z(s) \arrow[d, "i"'] \arrow[r, "i"] & X \arrow[d, "s"] \\
X \arrow[r, "0"]                    & E               
\end{tikzcd}\]
where $0$ is the zero section. The \textbf{localized Euler class} is defined in \cite[\S14.1]{fulton} as the class \begin{equation}
\label{loctopchern}
    e(E,s):=0^![X] \in A_{n-e}(Z(s)).
\end{equation}
It satisfies
\[i_*(e(E,s)) = c_e(E) \cap [X] \in A_{n-e}X . \]
For later use, we give the following proposition:
\begin{Pro}[Multiplicativity of localized Euler classes, {\cite[Example 17.4.8]{fulton}}]
\label{multiplicative}
    Let $E_i$ be a rank $e_i$ vector bundle over a pure $n$-dimensional scheme $X$ with a section $s_i$, where $i=1,2$, and $E=E_1 \oplus E_2,~s=s_1\oplus s_2$. Then,
    \[e(E,s) = e(E_1,s_1) \cup e(E_2,s_2) \in A_{n-(e_1+e_2)}(Z(s) ). \]
\end{Pro}

\subsection{Global construction of $Q_{g,m|k}(X_\Sigma,\beta)$}\label{globalconstruct} For a smooth project toric variety $X_\Sigma$ with its fan $\Sigma$, the quasimap moduli spaces $Q_{g,m|k}(X_\Sigma,\beta)$ can be embedded into a smooth Delign-Mumford stack over which there is a vector bundle with a section that singles out the quasimap moduli. In this way, one may construct the virtual fundamental class of the quasimap moduli space using the localized Euler class.

There are natural forgetful maps from $Q_{g,m|k}(X_\Sigma,\beta)$ to the $r$-fold fibered product of the Picard stack $\frak{Pic}_{g,m|k}^{r}$ with $r:=\rank\text{Pic}(X_\Sigma)$ by forgetting sections, and to the moduli stack of prestable curves $\frak{M}_{g,m|k}$ by forgetting both sections and line bundles. Denoted the forgetful maps by $\text{fgt}_{\text{s}}$ and $\text{fgt}_{\text{s,l}}$, respectively.
\[\begin{tikzcd}
Q_{g,m|k}(X_\Sigma,\beta) \arrow[r, "\text{fgt}_{\text{s}}"] \arrow[rd, "\text{fgt}_{\text{s,l}}"'] & \frak{Pic}_{g,m|k}^{r} \arrow[d, "\text{fgt}_{\text{l}}"] \\
                                          & \frak{M}_{g,m|k}                   
\end{tikzcd}\]
Write $\frak{Pic}^{r,\circ}$ and $\frak{M}^{\circ}$ as the substacks which the forgetful maps $\text{fgt}_{\text{s}}$ and $\text{fgt}_{\text{s,l}}$ factor through, respectively. Equivalently, $\frak{Pic}^{r,\circ}$ can be defined as the substack in the following way: For a family of quasimaps $C\rightarrow S$, there exists a positive integer $M=M(g,\beta,m,k)$ such that for all geometric points $s\in S$, $H^1(C_s,L_{\rho}(M) )=0$ for all $\rho$ \cite[Cor 3.1.11]{toricquasi}. Suppose that $\{\mathcal{P}_i\}$ is an integral basis for $\text{Pic}(X_\Sigma)$. We obtain the substack by imposing the following conditions:
\begin{enumerate}
    \item the degree of $L_i$ is equal to $\int_\beta c_1(\mathcal{P}_i)$
    \item the stability condition holds,
    \[\omega_C (\sum_{l=1}^{m}x_l+\sum_{l'=1}^{k}\epsilon y_{l'} ) \otimes \big( \otimes_{\rho \in \Sigma(1)} L_\rho^{\otimes \alpha_\rho} \big)^{\epsilon}\]
    is ample for every rational number $\epsilon > 0$, where $L_\rho:=\otimes_i L_i^{\otimes a_{i\rho}}$.
    \item $H^1(C,L_\rho(M))=0$ for all $\rho \in \Sigma(1)$.
\end{enumerate}  
These are open conditions, so the substack $\frak{Pic}^{r,\circ}$ is open.

Let $\pi_c : \mathcal{C} \rightarrow \frak{M}^{\circ}$ be the universal curve. Then, a universal curve, denoted by $\mathcal{C}\frak{Pic}^{\circ}$, over $\frak{Pic}^{r,\circ}$ is given by the fiber product of $\mathcal{C}$ and $\frak{Pic}^{r,\circ}$ over $\frak{M}^{\circ}$ with the two projections $\pi_1: \mathcal{C}\frak{Pic}^{r,\circ} 
\rightarrow \mathcal{C}$ and $\pi_2: \mathcal{C}\frak{Pic}^{r,\circ} 
\rightarrow \frak{Pic}^{r,\circ}$.
\[
\begin{tikzcd}
\mathcal{C}\frak{Pic}^{\circ} \arrow[r, "\pi_2"] \arrow[d, "\pi_1"'] & \frak{Pic}^{r,\circ} \arrow[d, "\pi_p"] \\
\mathcal{C} \arrow[r, "\pi_c"]                     & \frak{M}^{\circ}                   
\end{tikzcd}
\]
 By the stability condition (2) in Definition \ref{defquasimap} of a stable toric quasimap, one can take a $\pi_2$-relative amble bundle $\calO(1)$ over $\mathcal{C}\frak{Pic}^{\circ}$, whose fiber at $((C,\underline{x},\underline{y}),(L_i )_{i=1}^{r}) \in \frak{Pic}^{r,\circ}$ is given by 
\begin{equation}
\label{logcan}
    \omega_C (\sum_{l=1}^{m}x_l+\sum_{l'=1}^{k}\epsilon y_{l'} ) \otimes \big( \otimes_{\rho \in \Sigma(1)} L_\rho^{\otimes \alpha_\rho} \big)^{\epsilon_0},
\end{equation}
where $\epsilon_0$ can be chosen so that for all $\epsilon \in (0,\epsilon_0)\cap \Q$, such a line bundle with $\epsilon_0$ replaced by $\epsilon$ is isomorphic to the line bundle \eqref{logcan} from \cite[Cor 3.1.5 and Lem 3.1.10]{toricquasi}. Also, there are universal line bundles $\mathcal{L}_\rho$ over $\mathcal{C}\frak{Pic}^{\circ}$ with $\rho=1,\ldots, r$, whose fiber over a point $((C,\underline{x},\underline{y}),(L_i )_{i=1}^{r}) \in \frak{Pic}^{r,\circ}$ is $L_\rho$. Write $\mathcal{V}:=\oplus_\rho \mathcal{L}_\rho$.

As in the proof of \cite[Proposition 5]{Beh} or \cite[Lemma 2.5]{mark}, we take
\begin{equation}
\label{emb}
    \mathcal{B} :=  (\pi^{*}(\pi_{*}(\mathcal{V}^{\vee}(M) ) ) )^{\vee}(M),
\end{equation}
where $\pi:=\pi_2$ and $\mathcal{V}(M):=\mathcal{V}\otimes \calO (M)$. This gives an embedding \begin{equation}
\label{exactseq}
\mathcal{V} \hookrightarrow \mathcal{B} \twoheadrightarrow \mathcal{E}
\end{equation}
with $\mathbf{R}^1\pi_*( \mathcal{B} ) = 0$, where $\mathcal{E}$ is the cokernel of the embedding. Then, we also have $\mathbf{R}^1\pi_*( \mathcal{E} ) = 0$. Thus, $[\mathbf{R}\pi_*( \mathcal{B} ) 
\rightarrow \mathbf{R}\pi_*( \mathcal{E} )]$ forms a two-term resolution of $\mathbf{R}\pi_*( \mathcal{V} )$ by vector bundles.

Define the total space of sections of $\pi_*(\calB)$ by
\begin{equation}
\label{deftot}
\tot(\pi_*(\calB)):= \Spec (\Sym (\mathbf{R}^1 \pi_* (\omega_{\pi} \otimes \calB^{\vee}) ),
\end{equation}
where $\omega_\pi$ is the relative dualizing sheaf for $\pi$ and denote the corresponding map by \[p:\tot(\pi_*(\calB)) \rightarrow \frak{Pic}^{r,\circ}.\] The fiber of $p$ at $((C,\underline{x},\underline{y}),(L_i )_{i=1}^{r})$ is $\oplus_{\rho}H^0(C,L_\rho)$. Impose the generic nondegeneracy condition
appearing in Definition \ref{defquasimap}. Then, we obtain an open substack in $\tot(\pi_*(\calB))$ which is smooth and Deligne-Mumford, say $B$. By pulling back $\pi_*(\calE)$ via $p$ and restricting to $B$, we obtain a vector bundle $E:=p^* \pi_*(\calE)\lvert_{B}$, and a tautological section $s$ induced from the map
\[ H^0(\calC,\calB) \longrightarrow H^0(\calC,\calE). \]
The zero locus $Z(s)$ of the section $s$ is exactly the quasimap moduli space $Q_{g,m|k}(X_\Sigma,\beta)$.

This global model for a construction of the moduli space of quasimaps
\begin{equation}
\label{globalQ}
\begin{tikzcd}
               & E \arrow[d]                                  \\
Z(s) \arrow[r] & B \arrow[u, "s"', bend right, shift right=2]
\end{tikzcd}
\end{equation}
gives rise to the localized Euler class $e(E,s)$ defined in Section \ref{loceuler}. In \cite[Thm 3.2.1]{toricquasi}, the virtual fundamental class $[Q_{g,m|k}(X_\Sigma,\beta)]^{\text{vir}}$ is given in this way. This does not depend on the choice of embeddings and vector bundles and agrees with the Behrend--Fantechi virtual class \cite{Beh} defined using relative perfect obstruction $\mathbf{R}\pi_* \mathcal{B}^{\vee}$ (see \cite[Prop 2.14]{mark}). 

\section{The Quantum $H^*(T)$-module Structure on $H^*(X_\Sigma)$}
\label{sec_splitting}
\subsection{Quantum $H^*(T)$-module structure on $H^*(X_\Sigma)$}
Assume $g=0,~m=2$ and $k=1$ to define a $H^*(T)$-module structure on $H^*(X_\Sigma)$. Let $\{T_i\}$ be a basis for $H^*(X_\Sigma)$ and $\{T^j\}$ the dual basis under the intersection pairing on $X_\Sigma$.
\begin{Def}
\label{21oper}
For $\xi \in H^*(T)$ and $\phi \in H^*(X_\Sigma)$, define the \textbf{quantum $H^*(T)$-action} on $H^*(X_\Sigma)$ via the $2|1$-pointed quasimap invariants
\[\xi \star \phi := \sum_{\beta\in \text{Eff}}\sum_{i} \langle \phi, T_i \:|\xi \rangle_{0,2|1,\beta}T^i.\]
\end{Def}
The following theorem defines a $H^*(T)$-module structure on $H^*(X_\Sigma)$.
\begin{Thm}[Quantum $H^*(T)$-module structure]
\label{module}
For $\xi, \zeta \in H^*(T)$ and $\phi \in H^*(X_\Sigma)$,
\begin{equation}
\label{assoc}
\xi \star (\zeta \star \phi) = (\xi  \cdot \zeta ) \star \phi .
\end{equation}
\end{Thm}
This can be viewed as a modified version of WDVV equations via 3-pointed GW invariants. We need a splitting lemma to prove our $H^*(T)$-module structure in Theorem \ref{module}.

Expand the left-hand side of \eqref{assoc}:
\begin{align}
    \label{lhswdvv}
    \xi \star (\zeta \star \phi) &= \xi \star \big( \sum_{\beta\in \text{Eff}}\sum_{i}\langle \phi, T_i \:|\zeta \rangle_{0,2|1,\beta}T^i \big) \\ \nonumber
    &= \sum_{\beta \in \text{Eff}}\sum_{ \substack{\beta_1,\beta_2\in \text{Eff}\\\beta_1+ \beta_2=\beta} }\sum_{i,j} \langle \phi, T_i \:|\zeta \rangle_{0,2|1,\beta_1} \langle T^i, T_j \:|\xi \rangle_{0,2|1,\beta_2}T^j   
\end{align}
Similarly, the right-hand side of \eqref{assoc} is the following:
\begin{align}
    \label{rhswdvv}
    (\xi  \cdot \zeta ) \star \phi &=  \sum_{\beta\in \text{Eff}}\sum_{j}\langle \phi, T_j \:|\xi \cdot \zeta \rangle_{0,2|1,\beta}T^j
\end{align}
To have the equality between \eqref{lhswdvv} and \eqref{rhswdvv}, it is enough to see the degree $\beta$ part of the coefficients of $T^j$ are the same, i.e.,
\begin{equation}
\label{wdvv}
\sum_{ \substack{\beta_1,\beta_2\in \text{Eff}\\\beta_1+ \beta_2=\beta} }\sum_{i} \langle \phi, T_i \:|\zeta \rangle_{0,2|1,\beta_1} \langle T^i, T_j \:|\xi \rangle_{0,2|1,\beta_2}
=\langle \phi, T_j \:|\xi \cdot \zeta \rangle_{0,2|1,\beta}.
\end{equation}
\vspace{0.3cm}

There is a forgetful morphism \[\text{ft}:Q_{0,2|2}(X_\Sigma,\beta) \rightarrow \overline{M}_{0,2|2},\] where $\overline{M}_{0,2|d}$ is a Losev-Manin space defined in \cite{Losev}. Note that $\overline{M}_{0,2|2} \simeq \PP^1$, so that the divisors $D(13|24)$ and $D(3=4)$ are equivalent in $H^*(\overline{M}_{0,2|2})$, where
\begin{itemize}
    \item $D(13|24)$ is the class given by the locus of nodal curves that the one component has a heavy point marked by 1 and a light point marked by 3, and the other component has a heavy point marked by 2 and a light point marked by 4,
    \item $D(3=4)$ is the class given by the locus of an irreducible component where the 3rd and 4th light markings are colliding.
\end{itemize}
By pulling back through the forgetful morphism as in \cite[Prop 6.2.2]{abramo},
\begin{equation}
\label{forget}
    \text{ft}^{*}(D(13|24))=\text{ft}^{*}(D(3=4)).
\end{equation}
We will derive $\eqref{wdvv}$ from \eqref{forget} by showing that the left-hand side of \eqref{wdvv} is the same as 
\[ev_{1}^{*}(\phi)  ev_{2}^{*}(\psi) \hat{ev}_{3}^{*}(\xi) \hat{ev}_{4}^{*}(\zeta) \cap \text{ft}^{*}(D(13|24)),\]
and the right-hand side of \eqref{wdvv} is the same as 
\[ev_{1}^{*}(\phi)  ev_{2}^{*}(\psi) \hat{ev}_{3}^{*}(\xi) \hat{ev}_{4}^{*}(\zeta) \cap \text{ft}^{*}(D(3=4)).\]
The former part requires the diagonal pullback.

\subsection{Diagonal pullback} Let $\Delta: X_\Sigma \rightarrow X_\Sigma \times X_\Sigma$ be the diagonal embedding for a smooth projective toric variety $X_\Sigma$. It is regular of codimension equal to $\dim X_\Sigma$. For a basis $\{T_i\}$ and the dual basis $\{T^i\}$, \begin{equation}
\label{diagonal}[\Delta(X_\Sigma)] = \sum_{i} T_{i}\otimes T^i \in H^*(X_\Sigma)\otimes H^*(X_\Sigma) \simeq H^*(X_\Sigma \times X_\Sigma), \end{equation}
where $i$ in $T_i$ and $T^i$ means dimension and codimension, respectively. 

Recall the \eqref{globalQ} gives a global model for a construction of the virtual fundamental class of $Q:=Q_{0,2|2}(X_\Sigma,\beta)$. There is a gluing morphism
\begin{align*}
    \text{gl}:\mathfrak{M}_{0,2|1}^{\circ} \times \mathfrak{M}_{0,2|1}^{\circ} &\rightarrow \mathfrak{M}_{0,2|2}^{\circ} \\(C_1,(x_1,n_1;y_1))\times(C_2,(x_2,n_2;y_2)) &\mapsto (C_1\sqcup C_2 / n_1 \sim n_2,(x_1,x_2;y_1,y_2)),
\end{align*}
where $C_1\sqcup C_2 / n_1 \sim n_2$ is the nodal curve given by $C_1$ and $C_2$ gluing $n_1$ and $n_2$. Denote $Q_i:=Q_{0,2|1}(X_\Sigma,\beta_i),~i=1,2$. One can see that
\[ D(13|24;\beta_1,\beta_2) = Q_1 \times_{X_\Sigma}Q_2, \]
and the following fiber square commutes as in \cite[Prop 5.2.2]{abramo},
\[\begin{tikzcd}
\bigsqcup_{\beta_1 +\beta_2 = \beta} Q_1\times_{X_\Sigma}Q_2 \arrow[d] \arrow[r] & Q \arrow[d] \\
\mathfrak{M}_{0,2|1}^{\circ} \times \mathfrak{M}_{0,2|1}^{\circ} \arrow[r, "\text{gl}"]                    & \mathfrak{M}_{0,2|2}^{\circ}          
\end{tikzcd}\]
The following proposition relates the diagonal pullbacks of $[Q_1]^{\text{vir}} \times [Q_2]^{\text{vir}}$ for all $(\beta_1,\beta_2)$ with $[Q]^{\text{vir}}$ pulled back via the gluing morphism.
\begin{Pro}[Diagonal pullback]
\label{splittingaxiom}
The pullback of $[Q]^{\text{vir}}$ via the gluing morphism is the same as the diagonal pullback of $[Q_1]^{\text{vir}} \times [Q_2]^{\text{vir}}$, i.e., 
\[
\text{gl}^{!}[Q]^{\text{vir}} = \sum_{\beta_1 + \beta_2 =\beta}\Delta^! ([Q_1]^{\text{vir}} \times [Q_2]^{\text{vir}}),
\]
where $\Delta: X_\Sigma \rightarrow X_\Sigma \times X_\Sigma $ is the diagonal embedding.
\end{Pro}
Fix $\beta_1$ and $\beta_2$ with $\beta_1 + \beta_2 =\beta$ and let $D:= D(13|24;\beta_1,\beta_2)$. Denote a natural inclusion by $i_{\beta_1|\beta_2}:D \rightarrow \bigsqcup_{\beta_1 +\beta_2 = \beta} Q_1\times_{X_\Sigma}Q_2$. To prove Proposition \ref{splittingaxiom}, it is enough to show the following lemma.
\begin{Lem}
\label{comparison} The following holds
    \[i_{\beta_1|\beta_2}^{!} \textup{gl}^{!}[Q]^{\text{vir}}= \Delta^! ([Q_1]^{\text{vir}}\times [Q_2]^{\text{vir}}).\]
\end{Lem}

We pull back universal curves over $\mathfrak{Pic}_{0,2|1,\beta_1}^{r,\circ}$ via the projection $p_i$ and denote them by $\mathcal{C}_i$:
\[\begin{tikzcd}
\mathcal{C}_i  \arrow[d] \arrow[r, "\overline{p}_i"] & \mathcal{C}\mathfrak{Pic}_{0,2|1,\beta_i}^{r,\circ} \arrow[d] \\
\mathfrak{Pic}_{0,2|1,\beta_1}^{r,\circ} \times \mathfrak{Pic}_{0,2|1,\beta_2}^{r,\circ} \arrow[r, "p_i"]                    & \mathfrak{Pic}_{0,2|1,\beta_i}^{r,\circ}
\end{tikzcd}\]
Then, there is a universal curve \[\mathcal{C}_1 \bigsqcup \mathcal{C}_2 \rightarrow \mathfrak{Pic}_{0,2|1,\beta_1}^{r,\circ} \times \mathfrak{Pic}_{0,2|1,\beta_2}^{r,\circ},\]
with the inclusions $c_i : \calC_i \rightarrow  \mathcal{C}_1 \bigsqcup \mathcal{C}_2$.
    
Consider that there is a restriction morphism
    \begin{align*}
    r:\text{gl}^*\mathfrak{Pic}_{0,2|2,\beta}^{r,\circ} \longrightarrow \mathfrak{Pic}_{0,2|1,\beta_1}^{r,\circ} \times \mathfrak{Pic}_{0,2|1,\beta_2}^{r,\circ},
    \end{align*}
given in the following way: over $(C_i,(x_i,n_i;y_i))_{i=1}^{2}\in \mathfrak{M}_{0,2|1}^{\circ} \times \mathfrak{M}_{0,2|1}^{\circ}$,
    \begin{align*}
    (C_1\sqcup C_2 / n_1 \sim n_2,\{L_j\}_{j=1}^{r}) \longmapsto ((C_1,\{L_j\lvert_{C_1}\}_{j=1}^{r}),(C_2,\{L_j\lvert_{C_2}\}_{j=1}^{r})).
    \end{align*}
Denote $\text{gl}^*\mathcal{C}\mathfrak{Pic}_{0,2|2,\beta}^{r,\circ}$ by $\mathcal{C}$. There is a morphism $\overline{r}:r^* (\calC_1 \bigsqcup \calC_2) \rightarrow \calC_1 \bigsqcup \calC_2$ induced by the restriction morphism $r$. Also, there is a natural morphism
    \begin{equation*}
        \nu: r^{*}\big(\mathcal{C}_1 \bigsqcup \mathcal{C}_2 \big) \longrightarrow \mathcal{C},
    \end{equation*}
given by normalizing a nodal curve at the node. By abuse of notation, we write \[\nu_*\mathcal{V}_i:=\nu_*\overline{r}^*c_{i*}\overline{p}_i^*\mathcal{V}_i, \quad \nu_*\mathcal{B}_i:=\nu_*\overline{r}^*c_{i*}\overline{p}_i^*\mathcal{B}_i.\]

Observe that we have the normalization exact sequence
\[
\begin{tikzcd}
\text{gl}^*\mathcal{V} \arrow[r, hook]&
\bigoplus_{i=1}^{2}\nu_*\mathcal{V}_i \arrow[r, "d", two heads]&
\text{gl}^*\mathcal{V}\lvert_{n},
\end{tikzcd}
\]
where the map $d$ is given by the difference.

To have evaluation maps, using Lemma 3.4.1 in \cite{glsm}, we replace $\bigoplus_{i=1}^{2}\mathcal{B}_i $ by $\bigoplus_{i=1}^{2}\Tilde{\mathcal{B}}_i$ with a surjective morphism $\delta$ commuting the following diagram:
\begin{equation*}
    \begin{tikzcd}
\bigoplus_{i=1}^{2}\nu_* \mathcal{V}_i \arrow[d, "d"', two heads] \arrow[r, hook] & \bigoplus_{i=1}^{2}\nu_* \Tilde{\mathcal{B}}_i 
\arrow[ld, "\delta", two heads] 
\\
\text{gl}^*\mathcal{V}\lvert_{n}                                                               & 
\end{tikzcd}
\end{equation*}
with $\mathbf{R}^{1}\pi_*\Tilde{\mathcal{B}}_i=0$, where $\pi: \calC \rightarrow \text{gl}^*\mathfrak{Pic}_{0,2|2,\beta}^{r,\circ}$. Namely, define $\Tilde{\mathcal{B}}_i$ to be the equalizer of the following two morphisms:
\begin{align*}
&\mathcal{B}_i \oplus (\mathcal{V}_i)\lvert_{n_i} \xrightarrow{\text{proj}_1} \mathcal{B}_i    \xrightarrow{\text{restriction}} (\mathcal{B}_i )\lvert_{n_i}\\
&\mathcal{B}_i \oplus (\mathcal{V}_i)\lvert_{n_i} \xrightarrow{\text{proj}_2} (\mathcal{V}_i)\lvert_{n_i}    \xrightarrow{\text{inclusion}} (\mathcal{B}_i )\lvert_{n_i}
\end{align*}
Then, we obtain maps 
\begin{equation}
\label{evali}
    e_i:\Tilde{\mathcal{B}}_i\rightarrow  (\mathcal{V}_i)\lvert_{n_i}.
\end{equation}
Taking the difference of $e_1$ and $e_2$ defines the desired morphism $\delta$. Therefore, over $\mathcal{C}$, there is a commutative diagram
    \begin{equation}
        \label{commdiagbig}
        \begin{tikzcd}
\text{gl}^*\mathcal{V} \arrow[d, hook] \arrow[r, hook]                                         & \bigoplus_{i=1}^{2}\nu_*\Tilde{\mathcal{B}}_i \arrow[d, equal] \arrow[r, "s", two heads] & \Tilde{\mathcal{E}} \arrow[d, two heads]                  \\
\bigoplus_{i=1}^{2}\nu_*\mathcal{V}_i \arrow[d, "d"', two heads] \arrow[r, hook] & \bigoplus_{i=1}^{2}\nu_*\Tilde{\mathcal{B}}_i 
\arrow[ld, "\delta", two heads] 
\arrow[r, "s_1\oplus s_2", two heads]   & \bigoplus_{i=1}^{2}\nu_*\Tilde{\mathcal{E}}_i \\
\text{gl}^*\mathcal{V}\lvert_{n}                                                               &                                                                                                       &       
\end{tikzcd}
    \end{equation}
where $\Tilde{\mathcal{E}}$ and $\bigoplus_{i=1}^{2}\nu_*\Tilde{\mathcal{E}}_i$ are defined as the cokernels. 

The exact sequence in the second row defines $e(E_1\oplus E_2, s_1\oplus s_2)=[Q_1]^{\text{vir}} \times [Q_2]^{\text{vir}}$, where $E_1\oplus E_2$ is the pullback of $\pi_*\bigoplus_{i=1}^{2}\nu_*\Tilde{\mathcal{E}}_i$ along $p:\text{tot}(\pi_*\bigoplus_{i=1}^{2}\nu_*\Tilde{\mathcal{B}}_i) \rightarrow \text{gl}^*\mathfrak{Pic}_{0,2|2,\beta}^{r,\circ}$ as in the global construction method in Section \ref{globalconstruct}. The exact sequence in the first row in \eqref{commdiagbig} gives rises to $e(E,s)=\text{gl}^![Q]^{\text{vir}}$, where $E:=p^*\pi_*\Tilde{\mathcal{E}}$. The exact sequence
\begin{equation}
    \label{exactseqzdelta}
    \ker(\delta) \hookrightarrow \pi_*\bigoplus_{i=1}^{2}\nu_*\Tilde{\mathcal{B}}_i \twoheadrightarrow \text{gl}^*\mathcal{V}
\end{equation}
defines $e(V_n,\delta)$ where $V_n :=p^* \pi_*\text{gl}^*\mathcal{V} \lvert_n$. We give a relation among those three localized Euler classes.
\begin{Lem}
    \label{relatethree}
    The three localized Euler classes $e(E, s),~e(E_1\oplus E_2, s_1\oplus s_2),~e(V_n,\delta)$ are related in the following way:
    \[
    e(E, s)=e(E_1\oplus E_2, s_1\oplus s_2) \cup e(V_n,\delta) \in A^*(Z(s)).
    \]
\end{Lem}
 Before proving Lemma \ref{relatethree}, we give a lemma that will be used.
\begin{Lem}
\label{diagchase}
Given a commutative
\begin{equation*}
\begin{tikzcd}
                                                  &                                          &                                                                    & P  \\
A \arrow[rr, "\iota", hook, bend left=20] \arrow[rd, "\rho", hook] &                                          & M \arrow[r, "q", two heads] \arrow[rd, "\nu_{M}", two heads] \arrow[ru, "q_{P}", two heads] & B  &\quad  \\
                                                  & N \arrow[ru, "\iota_N", hook] \arrow[rr, "\nu_N", two heads] &                                                                    & P'
\end{tikzcd}
\end{equation*}
with exact sequences $A\hookrightarrow M \twoheadrightarrow B$, $N\hookrightarrow M \twoheadrightarrow P$ and $A \hookrightarrow N \twoheadrightarrow P'$, we can construct the following commutative
\begin{equation*}
\begin{tikzcd}[row sep=large,column sep=huge]
                                                                                                      & K' \arrow[rd, "\rho_{M}", hook]                               &                                                                                                      & P                                                                              &                                                        \\
A \arrow[rr, hook, bend left=20] \arrow[rd, hook] \arrow[ru, hook] \arrow[r, "\varphi"', "\sim", two heads, hook] & K \arrow[r, "\overline{\iota}"', hook] \arrow[d, "i_{N}", hook] \arrow[u, "i_{K'}", pos=0.7, hook] & M \arrow[r, two heads] \arrow[rd, two heads] \arrow[ru, two heads] \arrow[rr, "\overline{q}\quad", two heads, bend right=20] & B \arrow[u,"\Tilde{q}_P", two heads] \arrow[r, "\psi", "\sim"', two heads, hook] \arrow[d, "\nu_B", pos=0.7, two heads] & P\oplus P' \arrow[lu, "\pi_P"', two heads] \arrow[ld, "\pi_{P'}", two heads] \\
                                                                                                      & N \arrow[ru, hook] \arrow[rr, two heads]          &                                                                                                      & P'                                                                             &                                                       
\end{tikzcd}
\end{equation*}
with exact sequences $K \hookrightarrow M \twoheadrightarrow P\oplus P'$ and $K' \hookrightarrow M \twoheadrightarrow P'$.
\end{Lem}
\begin{proof}
    (Construction of $\Tilde{q}_{P}$ and its surjectivity) For $b\in B$, take $m\in q^{-1}(b)$, and define $\Tilde{q}_{P}(b):=q_{P}(m)$. It is well-defined; for some $m'\in q^{-1}(b)$, $q(m-m')=0$, thus $\iota(a)= m-m'$ for some $a \in A$. Since $\iota(a)=\iota_N ( \rho (a))$, By exactness of $N \hookrightarrow M \twoheadrightarrow P$, $q_P(m-m')=0$. Surjectivity follows by the construction.

    (Existence of an exact sequence $K \hookrightarrow M \rightarrow P\oplus P'$) By the universal property of product, there is a unique morphism $\overline{q}: M\rightarrow P\oplus P'$. Take $K:=\ker \overline{q}$.

    (Construction of $\psi: B \rightarrow P\oplus P'$ and $\nu_B$) Assume $a\in A$ with $q(\iota(a))=0$. Commutativity $\iota = \iota_N \circ \rho$ and exactness of $A \hookrightarrow N \twoheadrightarrow P'$, $q_P(\iota(a))=0$ and $\nu_M(\iota(a))=0$, respectively. Since $B=\coker \iota$, by the universal property of the cokernel, there is a unique morphism $\psi:B \rightarrow P\oplus P'$.

    (Surjectivity of $\nu_B$ and $\Tilde{q}_P$) Define $\nu_B:=\pi_{P'} \circ\phi$. $\nu_M = \nu_B \circ q$ and $q_P = \Tilde{q}_B \circ q$ imply that $\nu_B$ and $\Tilde{q}_P$ are surjective.

    (Construction of $K'$, $i_N$ and $i_{K'}$) Define $K':=\ker \nu_M$. For $k\in K$ with $\overline{q}\circ\overline{\iota}(k)=0$, then $q_P(\iota(k))=0$ and $\nu_M(\iota(k))=0$ by projecting via $\pi_P$ and $\pi_{P'}$. By the universal property of $N$ and $K'$ being the kernels, there are unique morphisms $i_N$ and $i_{K'}$. Since $\overline{\iota}$ is injective, one can show that $i_N$ and $i_{K'}$ are injective as well.

    (Construction of an isomorphism $\varphi$) We use the universal property of the kernel $K$. Note that $\nu_M(\iota(a))=0$ because of commutativity $\iota = \iota_N \circ \rho$ and $\nu_N = \nu_M\circ \iota_N$, and exactness of $A \hookrightarrow N \twoheadrightarrow P'$. Also, one can see that $q_P(\iota(a))=0$. From the construction of $\overline{q}$, $\overline{q}(\iota(a))=0$. Thus, there exist a unique map $\varphi:A \rightarrow K$ by the universal property. Next, we show that there is a unique morphism $\Tilde{\varphi}:K \rightarrow A$ such that $i_N = \rho \circ \Tilde{\varphi}$. We show that the composition $K \hookrightarrow N \twoheadrightarrow P'$ is the zero morphism. Consider that the composition $K \hookrightarrow M \twoheadrightarrow P'$ is the zero morphism since $K' \hookrightarrow M \twoheadrightarrow P'$ is exact and $\overline{\iota} = \rho_M \circ i_{K'}$. By commutativity $\iota_N \circ i_N = \overline{\iota}$ and $\nu_N = \nu_M\circ \iota_N$, we have the desired result. Therefore, $A$ being the kernel of $\nu_N$ gives rise to a unique morphism $\Tilde{\varphi}:K \rightarrow A$. Last, we show that $\varphi$ is an isomorphism. Consider that injectivity of $\rho$ implies injectivity of $\varphi$. For surjective, for $k \in K$, set $a:= \Tilde{\varphi}(k)$.
    Observe that \[i_N(\varphi(a)-k) = \rho(a)-i_N(k) = \rho(\Tilde{\varphi}(k))-i_N(k) = i_N(k)-i_N(k)=0.\]
    Since $i_N$ is injective, $\varphi(a)=k$, which implies that $\varphi$ is surjective.

    (Claim: $\psi$ is an isomorphism) Since $A \simeq K$, $B\simeq \text{Im}(\overline{q})$. Thus, it is enough to show that $\psi$ is surjective. Let $(p,p') \in P \oplus P'$. Using $\pi_P$ and $\pi_{P'}$, there exist $m \in M$ and $n \in N$ such that $\overline{q}(m)=(p,x),~\overline{q}(\iota_N (n)) =(y,p')$ where $x:=\nu_M(m)$ and $y:=q_P(\iota_N (n))$. Consider that exactness of $N \hookrightarrow M \twoheadrightarrow P$ implies that $y=0$. The same way shows that there exists $n' \in N$ such that $\overline{q}(\iota_N (n')) =(0,x)$. Then, observe that
    \[\overline{q}(m+\iota_N (n) - \iota_N (n') ) = (p,x) + (0,p') - (0,x) = (p,p'). \]
    Hence, $\psi$ is an isomorphism.
\end{proof}
We apply Lemma \ref{diagchase} to the diagram \eqref{commdiagbig}.
\begin{proof}[Proof of Lemma \ref{relatethree}]
    Take
    \begin{align*}
        A&:=\text{gl}^*\mathcal{V},&M&:=\bigoplus_{i=1}^{2}\nu_*\Tilde{\mathcal{B}}_i,&B&:=\Tilde{\mathcal{E}}\\
        N&:=\bigoplus_{i=1}^{2}\nu_*\mathcal{V}_i,&P&:=\bigoplus_{i=1}^{2}\nu_*\Tilde{\mathcal{E}}_i,&P'&:=\text{gl}^*\mathcal{V}\lvert_{n}
    \end{align*}
    By Lemma \ref{diagchase},
\begin{equation*}
    \Tilde{\calE}\simeq \bigg(\bigoplus_{i=1}^{2}\nu_*\Tilde{\mathcal{E}}_i \bigg) \bigoplus \text{gl}^*\mathcal{V}\lvert_{n}.
\end{equation*}
Apply Proposition \ref{multiplicative} to obtain the relation among three classes.
\end{proof}
Recall that $e(E, s)=[Q]^{\text{vir}}$, $e(E_1\oplus E_2, s_1\oplus s_2)=[Q_1]^{\text{vir}}\times [Q_2]^{\text{vir}}$. We prove Lemma \ref{comparison}
\begin{proof}[Proof of Lemma \ref{comparison}] Let $Y:=\text{tot}\big(\bigoplus_{i=1}^{2}\nu_* \Tilde{\mathcal{B}}_i \big)$. Note that $Y$ is smooth.

(Step1: $\Delta^!{[Y]} = e(V_n,\delta)$) Consider the following diagram of fiber squares:
    \[\begin{tikzcd}
Z(\delta) \arrow[r,"\text{id}"] \arrow[d,"\text{id}"'] & Z(\delta) \arrow[r] \arrow[d] & Y \arrow[d, "0_{V_n}"] \\
Z(\delta) \arrow[r] \arrow[d] & Y \arrow[r,"\delta"] \arrow[d] & V_n \\
X_\Sigma \arrow[r,"\Delta"] & X_\Sigma \times X_\Sigma &  
\end{tikzcd}
\]
The evaluation map $Y \rightarrow [V/T] \times [V/T]$ at $n_1$ and $n_2$ is induced from $e_1$ and $e_2$ in \eqref{evali}. The evaluation map $Z(\delta)\rightarrow [V/T]$ is well-defined since $Z(\delta)$ is defined where the evaluation maps induced from $e_1$ and $e_2$ coincide. Since $n_1$ and $n_2$ do not collide with any base points and light markings, those evaluation maps factor through $X_\Sigma \times X_\Sigma$ and $X_\Sigma$, respectively. Then, from \cite[Theorem 6.4]{fulton},
\begin{equation}
\label{pullbackT1}
    0_{V_n}^!\Delta^! [Y] = \Delta^!0_{V_n}^! [Y]
\end{equation}
Using the definition of the localized top Chern class \eqref{loctopchern} and \cite[Remark 6.2.1]{fulton}, the equation \eqref{pullbackT1} becomes
\begin{align*}
     \text{id}^*\Delta^! [Y]&= \text{id}^*e(V_n,\delta)\\
     \Delta^! [Y]&= e(V_n,\delta)
\end{align*}
(Step2: $\Delta^!e(E_1\oplus E_2, s_1 \oplus s_2) = e(E_1\oplus E_2, s_1 \oplus s_2) \cup e(V_n,\delta)$) Consider the following diagram of fiber squares:
\[\begin{tikzcd}
Z(s) \arrow[r] \arrow[d] & Z(s_1\oplus s_2) \arrow[r] \arrow[d] & Y \arrow[d, "0_{E_1 \oplus E_2}"] \\
Z(\delta) \arrow[r] \arrow[d] & Y \arrow[r,"s_1 \oplus s_2"] \arrow[d] & E_1 \oplus E_2 \\
X_\Sigma \arrow[r,"\Delta"] & X_\Sigma \times X_\Sigma &  
\end{tikzcd}
\]
We apply \cite[Theorem 6.4]{fulton}, again.
\begin{equation}
\label{pullbackT2}
    0_{E_1 \oplus E_2}^!\Delta^! [Y] = \Delta^!0_{E_1 \oplus E_2}^! [Y] = \Delta^!e(E_1\oplus E_2, s_1 \oplus s_2)
\end{equation}
By the Step1, the left-hand side of \eqref{pullbackT2} becomes $0_{E_1 \oplus E_2}^!e(V_n,\delta)$. Note that surjectivity of $\delta: \oplus\Tilde{\mathcal{B}_i} \rightarrow \text{gl}^*\mathcal{V}\lvert_n$ implies that $Z(\delta)$ is smooth. 
Thus, 
\[
0_{E_1 \oplus E_2}^!e(V_n,\delta) = e(E_1\oplus E_2, s_1 \oplus s_2) \cup e(V_n,\delta).
\]

Last, applying Lemma \ref{relatethree}, we conclude with $\Delta^! e(E_1\oplus E_2, s_1 \oplus s_2) = e(E,s)$.
\end{proof}

\subsection{The pullback $\text{ft}^{*}(D(13|24))$} We describe splitting lemma for $\text{ft}^{*}(D(13|24))$. The class $\text{ft}^{*}(D(13|24))$ is given from \[\text{ft}^{-1}(D(13|24)) = \bigsqcup_{\beta_1 + \beta_2 =\beta}D(13|24;\beta_1|\beta_2),\]
where $D(13|24;\beta_1,\beta_2)$ is the locus of quasimaps whose source curve is a nodal curve where the first and third markings lie in some components whose degree is $\beta_1$, and the second and fourth markings lie in the other components of degree $\beta_2$.

We give a splitting lemma to prove the WDVV equations \eqref{wdvv}:
\begin{Lem}[Splitting lemma]
\label{split}
    For $\beta_1,\beta_2 \in \text{Eff}$ with $\beta_1 + \beta_2 = \beta$, $\phi,\psi \in H^{*}(X_\Sigma)$, and $\xi, \zeta \in H^{*}([V/T])$,
    \begin{equation}
    \label{spliteqn}
        \int_{[D(13|24;\beta_1,\beta_2)]} ev_{1}^{*}(\phi)  ev_{2}^{*}(\psi) \hat{ev}_{3}^{*}(\xi) \hat{ev}_{4}^{*}(\zeta) = \sum_{i} \langle \phi, T_i \:|\xi \rangle_{0,2|1,\beta_1} \langle \psi, T^i \:|\zeta \rangle_{0,2|1,\beta_2}
    \end{equation}
\end{Lem}
Our strategy to prove Lemma \ref{split} is based on \cite{abramo}. Set $X:=X_\Sigma$. One can write down
\[ev_{2*}(ev_1^*(\phi)\hat{ev}_3^*(\xi) \cap [Q_1]^{\text{vir}} ) = \sum_{j}c_jT^j.\]
Applying the projection formula gives us
\begin{align*}
    \langle \phi,T_i \mid \xi \rangle_{\beta_1} &= ev_{2*}(ev_1^*(\phi)ev_2^*(T_i)\hat{ev}_3^*(\xi) \cap [Q_1]^{\text{vir}} ) \cap X\\
    &=  T_i ev_{2*}(ev_1^*(\phi)\hat{ev}_3^*(\xi) \cap [Q_1]^{\text{vir}} ) \cap X=c_i.
\end{align*}
Write $\overline{ev}$ for the evaluations on $Q_2$. Then, it is easy to see that
\begin{align}
    &\sum_{i}\langle \phi,T_i \mid \xi \rangle_{\beta_1} \langle \psi,  T^i \mid  \zeta \rangle_{\beta_2} =\langle \psi,  \sum_{i}\langle \phi,T_i \mid \xi \rangle_{\beta_1}T^i \mid  \zeta \rangle_{\beta_2} \nonumber \\ \label{evlower}
    &=  \overline{ev}_1^*(\psi) \overline{\hat{ev}}_3^*(\zeta) \overline{ev}_2^* ( ev_{2*}(ev_1^*(\phi)\hat{ev}_3^*(\xi) \cap [Q_1]^{\text{vir}} ) )   \cap [Q_2]^{\text{vir}}.
\end{align}
Using the abbreviation $e$ for $ev$, consider the following diagram:
\[
\begin{tikzcd}
Q_1 \times_{X} Q_2 \arrow[d,"p_1"'] \arrow[r, "p_2"]              & Q_2 \arrow[d, "\overline{e}_2"] \\
Q_1 \arrow[r, "e_2"]                                         & X                         
\end{tikzcd}
\]
We prove the following lemma, which is Lemma 6.2.7 in \cite{abramo} modified to our case.
\begin{Lem}
    For $z_i \in A^*(Q_i)$, the following equality holds:
    \begin{equation}
    \label{eqabra}
        p_{2*} ( p_1^* z_1  p_2^* z_2 \cap \Delta^!( [Q_1]^{\text{vir}} \times [Q_1]^{\text{vir}} )) =   (\overline{e}_2^*  e_{2*}(z_1 \cap [Q_1]^{\text{vir}} ) )  z_2 \cap [Q_2]^{\text{vir}}
    \end{equation}
\end{Lem}
\begin{proof}
For convenience, write $e:=e_2$ and $\overline{e}:=\overline{e}_2$. Consider the following diagram:
\[
\begin{tikzcd}
Q_1 \times_{X} Q_2 \arrow[d] \arrow[r, "p_2"]              & Q_2 \arrow[r, "\overline{e}"] \arrow[d, "\Gamma_{\overline{e}}"]               & X \arrow[d, "\Delta"] \\
Q_1 \times Q_2 \arrow[r, "e\times id"] \arrow[d, "\pi_1"'] & X \times Q_2 \arrow[r, "id \times \overline{e}"] \arrow[d, "\overline{\pi}_1"] & X \times X            \\
Q_1 \arrow[r, "e"]                                         & X                                                                              &                      
\end{tikzcd}
\]
where $\Gamma_{\overline{e}}$ is the graph of $\overline{e}$. Denote the projections on the second factor by $\pi_2:Q_1 \times Q_2 \rightarrow Q_2$ and $\overline{\pi}_2:X \times Q_2 \rightarrow Q_2$. Then, $\overline{\pi}_2 \circ\Gamma_{\overline{e}} =\id$ and $\overline{\pi}_1 \circ\Gamma_{\overline{e}} =\overline{e}$. Since $X$ is smooth, $\Gamma_{\overline{e}}$ is regular. Thus, one can apply \cite[Thm 6.2(c), Thm 6.2(a) and Remark 6.2.1]{fulton} and the projection formula to the left-hand side of \eqref{eqabra}. Then,
\begin{align*}
    p_{2*} ( p_1^* z_1  p_2^* z_2 \cap \Delta^!( [Q_1]^{\text{vir}} \times [Q_2]^{\text{vir}} )) &=  p_{2*} \Delta^!(  \pi_1^* z_1 \pi_2^* z_2 \cap ( [Q_1]^{\text{vir}} \times [Q_2]^{\text{vir}} )) \\
    &= p_{2*} \Gamma_{\overline{e}}^!( \pi_1^* z_1  \pi_2^* z_2 \cap ( [Q_1]^{\text{vir}} \times [Q_2]^{\text{vir}} ) )\\
    &=\Gamma_{\overline{e}}^! (e\times id)_{*} ( \pi_1^* z_1  \pi_2^* z_2 \cap ( [Q_1]^{\text{vir}} \times [Q_2]^{\text{vir}} )  )\\
    &=\Gamma_{\overline{e}}^* (e\times id)_{*} ( \pi_1^* z_1  \pi_2^* z_2 \cap ( [Q_1]^{\text{vir}} \times [Q_2]^{\text{vir}} ))\\
    &=\Gamma_{\overline{e}}^* (e\times id)_{*} ( \pi_1^* (z_1  \cap  [Q_1]^{\text{vir}}) \pi_2^* (z_2 \cap  [Q_2]^{\text{vir}})  )\\
    &=\Gamma_{\overline{e}}^* ( \overline{\pi}_1^* e_* (z_1  \cap  [Q_1]^{\text{vir}}) \overline{\pi}_2^* (z_2 \cap  [Q_2]^{\text{vir}})  )\\
    &=\Gamma_{\overline{e}}^*\overline{\pi}_1^* e_* (z_1  \cap  [Q_1]^{\text{vir}}) \Gamma_{\overline{e}}^*\overline{\pi}_2^* (z_2 \cap  [Q_2]^{\text{vir}})  )\\
    &=(\overline{e}^* e_* (z_1  \cap  [Q_1]^{\text{vir}}))z_2 \cap  [Q_2]^{\text{vir}})  
\end{align*}
\end{proof}

We prove Lemma \ref{split}.
\begin{proof}[proof of Lemma \ref{split}]
Take $z_1:=ev_1^*(\phi)\hat{ev}_3^*(\xi)$ and $z_2:=\overline{ev}_1^*(\psi) \overline{\hat{ev}}_3^*(\zeta)$ and apply Lemma \ref{eqabra}. Then, \eqref{evlower} becomes
\begin{equation}
\label{eqapplydiag}
p_{2*} ( p_1^* (ev_1^*(\phi)\hat{ev}_3^*(\xi))  p_2^* (\overline{ev}_1^*(\psi) \overline{\hat{ev}}_3^*(\zeta)) \cap \Delta^!( [Q_1]^{\text{vir}} \times [Q_2]^{\text{vir}} )).
\end{equation}
Recall that $D:= D(13|24;\beta_1,\beta_2)$, $Q:=Q_{0,2|2}(X_\Sigma,\beta)$, $Q_i:=Q_{0,2|1}(X_\Sigma,\beta_i)$ for $i=1,2$ and $i_{\beta_1|\beta_2}:D \rightarrow \bigsqcup_{\beta_1 +\beta_2 = \beta} Q_1\times_{X_\Sigma}Q_2$. Applying Lemma \ref{comparison} to \eqref{eqapplydiag} gives us
\begin{equation}
\label{eqapplydiag1}
p_{2*} ( p_1^* (ev_1^*(\phi)\hat{ev}_3^*(\xi))  p_2^* (\overline{ev}_1^*(\psi) \overline{\hat{ev}}_3^*(\zeta)) \cap i_{\beta_1|\beta_2}^{!} \textup{gl}^{!}[Q]^{\text{vir}} ).
\end{equation}
Pushing forward \eqref{eqapplydiag1} to the Chow ring of a point, we obtain \eqref{spliteqn} in Lemma \ref{split}. Therefore, we proved Lemma \ref{split}.
\end{proof}
Therefore, by Lemma \ref{split}, the left-hand side of \eqref{wdvv} becomes
\begin{align}
\label{ftd1324}
    &\sum_{  \substack{\beta_1,\beta_2\in \text{Eff}\\\beta_1+ \beta_2=\beta} 
 }\int_{[D(13|24;\beta_1,\beta_2)]} ev_{1}^{*}(\phi)  ev_{2}^{*}(\psi) \hat{ev}_{3}^{*}(\xi) \hat{ev}_{4}^{*}(\zeta) \nonumber \\ 
 &=\int_{\text{ft}^*(D(13|24))} ev_{1}^{*}(\phi)  ev_{2}^{*}(\psi) \hat{ev}_{3}^{*}(\xi) \hat{ev}_{4}^{*}(\zeta)
\end{align}

\subsection{The pullback $\text{ft}^{*}(D(3=4))$}  Consider that there is a natural inclusion \[\text{inc}:\frak{M}_{0,2|1}^{\circ} \hookrightarrow \frak{M}_{0,2|2}^{\circ},\]
by sending $(C;x_1,x_2,y_1)$ to $(C;x_1,x_2,y_1,y_1)$, since two light points in $\frak{M}_{0,2|2}^{\circ}$ can collide. Denote $Q':=Q_{0,2|1}(X_\Sigma,\beta)$. There is a fibered square
\[\begin{tikzcd}
Q' \arrow[d] \arrow[r] & Q \arrow[d] \\
\mathfrak{M}_{0,2|1}^{\circ} \arrow[r, "\text{inc}"]                    & \mathfrak{M}_{0,2|2}^{\circ}          
\end{tikzcd}\]
We successively pull back along the inclusion to obtain the following diagram:
\[\begin{tikzcd}
\mathcal{V} \arrow[d] & \arrow[l,"\text{inc}"'] \text{inc}^*\mathcal{V} \arrow[d]  \\
\mathcal{C}\frak{Pic} _{0,2|2}^{r,\circ} \arrow[d,"\pi"'] & \arrow[l,"\text{inc}"'] \text{inc}^*\mathcal{C}\frak{Pic} _{0,2|2}^{r,\circ} \arrow[d,"\pi'"'] \\
\frak{Pic} _{0,2|2}^{r,\circ} \arrow[d] & \arrow[l,"\text{inc}"'] \text{inc}^*\frak{Pic} _{0,2|2}^{r,\circ} \arrow[d] \\
\frak{M}_{0,2|2}^{\circ} & \arrow[l,"\text{inc}"'] \frak{M}_{0,2|1}^{\circ}
\end{tikzcd}\]
For the notations, we use inc for all the induced inclusions. One can observe that the embedding $\mathcal{V} \rightarrow \mathcal{B}$ pulled back over $\frak{M}_{0,2|1}^{\circ}$ via the inclusion defines the virtual fundamental classes $[Q']^{\text{vir}}$ using the method in Section \ref{globalconstruct}. By the definition of the total space of sections \eqref{deftot}, we have an induced map
\[Y':=\text{tot}(\pi'_{*}\mathcal{B'}) \rightarrow Y:=\text{tot}(\pi_{*}\mathcal{B}),\]
where $\mathcal{B}':=\text{inc}^*\mathcal{B} $. Also, denote $\mathcal{E}':=\text{inc}^*\mathcal{E} $ and $E':= p'^*\pi'_*(\mathcal{E}')$, where $p': \text{tot}(\pi_{*}^{*}\mathcal{B'}) \rightarrow  \text{inc}^*\frak{Pic} _{0,2|1}^{r,\circ}$. There are induced maps $Q' \rightarrow Q$ and $E' \rightarrow E$.

We relate the virtual fundamental classes $[Q]^{\text{vir}}$ and $[Q']^{\text{vir}}$.
\begin{Lem}
\label{inclusionpullback}
    The following holds
    \[\textup{inc}^![Q]^{\text{vir}} = [Q']^{\text{vir}}.\]
\end{Lem}
\begin{proof}
From \cite[Thm 6.4]{fulton}, the fibered diagram
    \[
\begin{tikzcd}
Q' \arrow[d] \arrow[r]              & Q \arrow[r] \arrow[d]               & Y \arrow[d, "0"] \\
Y' \arrow[r] \arrow[d] & Y  \arrow[r] \arrow[d] & E            \\
\frak{M}_{0,2|1}^{\circ} \arrow[r, "\text{inc}"]                                         & \frak{M}_{0,2|2}^{\circ}                                                                              &                      
\end{tikzcd}
\]
 implies
\[\textup{inc}^![Q]^{\text{vir}} = \textup{inc}^!0^![Y] = 0^!\textup{inc}^![Y].\]
Then, since the zero sections $0$ and $0'$ are regular embeddings, by \cite[Thm 6.2(c)]{fulton}, the following fibered squares
\[
\begin{tikzcd}
Q' \arrow[d] \arrow[r]              & Y' \arrow[d]               \\
Y' \arrow[r,"0'"] \arrow[d] & E'  \arrow[d]             \\
Y \arrow[r, "0"]                                         & E  
\end{tikzcd}
\]
 gives rise to
\[0^!\textup{inc}^![Y] = 0'^!\textup{inc}^![Y] = 0'^![Y'] = [Q']^{\text{vir}}.\]
\end{proof}

Consider the following commutative diagram
\[
\begin{tikzcd}[column sep=7em]
Q  \arrow[r,"ev_1\times ev_2\times \hat{ev}_3 \times \hat{ev}_4"]              & X_\Sigma \times X_\Sigma \times [V/T] \times [V/T] \\
Q' \arrow[u]\arrow[r,"ev'_1\times ev'_2\times \hat{ev}'_3 \times \hat{ev}'_3"]  & X_\Sigma \times X_\Sigma \times [V/T] \times [V/T]  \arrow[u,"id"]          
\end{tikzcd}
\]
Then, applying Lemma \ref{inclusionpullback} gives us
\begin{align}
\label{eqnincpullback}
&\text{inc}^! (ev_{1}^{*}(\phi)  ev_{2}^{*}(\psi) \hat{ev}_{3}^{*}(\xi) \hat{ev}_{4}^{*}(\zeta) \cap [Q]^{\text{vir}})    \\
&=\text{inc}^! (ev_{1}^{*}(\phi)  ev_{2}^{*}(\psi) \hat{ev}_{3}^{*}(\xi) \hat{ev}_{4}^{*}(\zeta) ) \cap \text{inc}^! [Q]^{\text{vir}} \nonumber\\
&={ev'}_{1}^{*}(\phi)  {ev'}_{2}^{*}(\psi) \hat{ev}{'}_{3}^{*}(\xi) \hat{ev}{'}_{3}^{*}(\zeta)  \cap [Q']^{\text{vir}} \nonumber\\
&={ev'}_{1}^{*}(\phi)  {ev'}_{2}^{*}(\psi) \hat{ev}{'}_{3}^{*}(\xi \cup \zeta)  \cap [Q']^{\text{vir}} \nonumber
\end{align}
Therefore, since the following fibered square 
\[\begin{tikzcd}
Q'  \arrow[r] \arrow[d]             & Q \arrow[d,"\text{ft}"] \\
D(3=4) \arrow[r]  & \overline{M}_{0,2|2}           
\end{tikzcd}\]
implies $\text{ft}^*(D(3=4))\simeq Q'$, by \eqref{eqnincpullback}, the right-hand side of \eqref{wdvv} becomes
\[\langle \phi, T_j \:|\xi \cdot \zeta \rangle_{0,2|1,\beta}= \int_{\text{ft}^{*}(D(3=4))}ev_{1}^{*}(\phi)  ev_{2}^{*}(\psi) \hat{ev}_{3}^{*}(\xi) \hat{ev}_{4}^{*}(\zeta).\]
Hence, from the equivalence \eqref{forget}, we achieve the equality in \eqref{wdvv}.

\section{The Quantum $H^*((\C^*)^2)$-module for the Hirzebruch Surface $\F_2$ }\label{sec_compute}
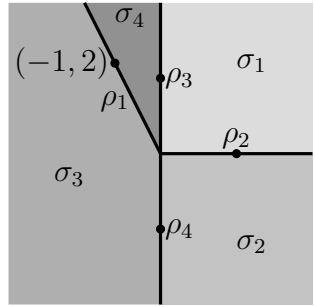
\begin{figure}[b] 
    \centering
\begin{tikzpicture}
\fill[color = {rgb,255:red,220; green,220; blue,220}] (0,0) -- (2,0) -- (2,2) --(0,2)--(0,0);
\fill[color = {rgb,255:red,145; green,145; blue,145}] (0,0) -- (0,2) -- (-1,2)--(0,0);
\fill[color = {rgb,255:red,175; green,175; blue,175}] (0,0) -- (-1,2) -- (-2,2) --(-2,-2)--(0,-2)--(0,0);
\fill[color = {rgb,255:red,195; green,195; blue,195}] (0,0) -- (2,0) -- (2,-2) --(0,-2)--(0,0);
\draw[very thick] (0,0) --(2,0);
\draw[very thick] (0,0) --(0,2);
\draw[very thick] (0,0) --(-1,2);
\draw[very thick] (0,0) --(0,-2);

\node at (1.2,1.2) {$\sigma_1$};
\node at (1.2,-1.2) {$\sigma_2$};
\node at (-1.2,-0.3) {$\sigma_3$};
\node at (-0.4,1.8) {$\sigma_4$};

\node at (-1.3,1.2) {$(-1,2)$};
\filldraw (-0.6,1.2) circle (1.5pt);
\node at (-0.6,0.7) {$\rho_1$};

\filldraw (0,1) circle (1.5pt);\node at (0.25,1) {$\rho_3$};

\filldraw (1,0) circle (1.5pt);\node at (1,0.2) {$\rho_2$};

\filldraw (0,-1) circle (1.5pt);\node at (0.25,-1) {$\rho_4$};
\end{tikzpicture}
\caption{The fan for $\F_2$}
\label{fanhirz}
\end{figure}

In this section, we compute the $H^*((\C^*)^2)$-module structure for the Hirzebruch surface of type 2, which we denote by $\F_2$. It may be represented as a toric variety with fan $\Sigma$ depicted in Figure \ref{fanhirz}, where $\rho_1=(-1,2),~ \rho_2=(1,0),~ \rho_3=(0,1),~ \rho_4=(0,-1)$ are primitive ray generators.
Note that $|\Sigma(1)|=4$ and $\rank\Pic(\F_2)=2$. Thus, there is a $(\C^*)^2$-action on $\C^4$ given by a $4 \times 2$ matrix; in this case, the matrix is \begin{equation}
\label{actionmatrix}
    \begin{pmatrix}
1&1&2&0\\
0&0&1&1
\end{pmatrix}.
\end{equation}  The geometric quotient description from Section \ref{sec_pre} for $\F_2$ is 
\[\F_2 \simeq \C^4\git (\C^*)^2.\]
Each column of the matrix gives rise to a line bundle $\calO_{\F_2}(D_i)$, where $D_i$ is the prime torus invariant divisor on $\F_2$ defined as the zero locus of the $i$-th coordinate function on $\C^4$. We omit the subscript $\F_2$ in the notation of such a line bundle when the context is clear.

For the complete description of the $H^*((\C^*)^2)$-module structure, it is required to compute all $2|1$-pointed quasimap invariants for $\F_2$. We will apply the Atiyah--Bott localization formula in \cite[\S4.3]{mirrsymm}.

\subsection{Linearization} From the theory of equivariant cohomology proposed in \cite[\S2]{fultonequiv}, the $\C^*$-equivariant cohomology of a point or the group cohomology of $\C^*$, denoted by $H^{*}(\C^*):=H_{\C^*}^{*}(\text{pt})$, is the polynomial ring $\C[\alpha]$ where $\alpha$ represents the Chern class of the hyperplane line bundle $\calO_{\PP^\infty}(1)$ (this convention is from \cite[\S27.1]{mirrsymm}). Consider the complex line bundle over a point $L_k \rightarrow \text{pt}$ with a $\C^*$-action on the fiber given by $t.v=t^k v$. With the choice of the convention for $\alpha$, the equivariant top Chern class of $L_k$, denoted by $e^{\C^*}(L_k)$, is $-k\alpha$. We call this the \textit{weight} of $L_k$.

Consider a $\mathcal{T}:=(\C^*)^4$-action on $\F_2\simeq \C^4\git (\C^*)^2$ given by \[(t_1,t_2,t_3,t_4).(z_1,z_2,z_3,z_4) = (t_1z_1,t_2z_2,t_3z_3,t_4z_4).\]
The corresponding fixed points are
\[p_1=(1:0:1:0),~p_2=(0:1:1:0),~p_3=(0:1:0:1),~p_4=(1:0:0:1).\]
One can find the weights of the tangent spaces at each $(\C^*)^4$-fixed point with respect to this action. For example, at the fixed point $p_1=(1:0:1:0) = D_2 \cap D_4$, we have $(1:\frac{t_2}{t_1}:1:\frac{t_1^2t_4}{t_3})$. It follows that
\[e^\mathcal{T}(\mathcal{T}_{p_1} \F_2)= e^\mathcal{T}(\calO^\mathcal{T} (D_2)\rvert_{p_1} \oplus \calO^\mathcal{T} (D_4)\rvert_{p_1}) = (\alpha_1-\alpha_2) \cdot (\alpha_3-\alpha_4-2\alpha_1),   \]
where $\alpha_1-\alpha_2$ and $\alpha_3-\alpha_4-2\alpha_1$ are the weights at $p_1$ along $D_4$-direction and $D_2$-direction, respectively, \cite[\S4.3]{melissa}. In the same fashion, the four line bundles $\calO (D_i)$ are canonically linearized in a sense that the weights of $\calO^\mathcal{T} (D_i)\lvert_{p_k}$ are given by the weights of the corresponding tangent spaces at the fixed point. All the weights of the equivariant line bundles associated with $D_i$ restricted to each fixed point are in table \ref{weighttable}.

From the fan $\Sigma$ of $\F_2$, the cohomology ring of $\F_2$ can be written using the divisors $D_2$ and $D_4$ as generators: \[H^*(\F_2) \simeq \Q[D_2,D_4]/\langle D_2^2,D_4^2+2D_2D_4 \rangle,\]
where $D_2D_4 =[\text{pt}]$. The reason why the divisors $D_2$ and $D_4$ were chosen is discussed in Section \ref{subsecbat}. The cohomology ring has a graded structure as a $\C$-vector space: \[H^*(\F_2) = H^0(\F_2) \oplus H^2(\F_2) \oplus H^4(\F_2)  = \C\cdot [\F_2] \oplus (\C\cdot [D_2] \oplus \C\cdot [D_4])  \oplus \C\cdot [pt].\]

The cohomology of the corresponding stack quotient is given as follows:
\[H^*([\C^4 / (\C^*)^2]) \simeq H^*((\C^*)^2) \simeq \C[\sigma_2,\sigma_4],\]
where $\sigma_2$ and $\sigma_4$ are the classes from the line bundles determined by the second column $(1,0)^T$ and fourth column $(0,1)^T$ of the action matrix \eqref{actionmatrix} above.

The following Theorem describes an explicit quantum $H^*((\C^*)^2)$-module structure on $H^*(\F_2)$ given from Theorem \ref{module}.
\begin{Thm}
\label{Thm}
    The quantum $H^*((\C^*)^2)$-module structure for $\F_2$ (as defined in a previous section) is given by the following
\begin{align*}
    \sigma_2 \star 1 &= D_2-\frac 12 f(q_4)D_4 &\quad \sigma_4 \star 1&= (1+f(q_4))D_4 \\
    \sigma_2 \star D_2 &=q_2q_4(1+f(q_4)) -\frac{1}{2}f(q_4)pt &\quad \sigma_4 \star D_2&=-\frac{1}{2}q_2f(q_4) + (1+f(q_4))pt\\
    \sigma_2 \star D_4 &=-2q_2q_4(1+f(q_4)) + (1+f(q_4))pt &\quad \sigma_4 \star D_4&=q_2(1 + f(q_4))-2(1+f(q_4))pt\\
    \sigma_2 \star pt &= q_2q_4(1+f(q_4))D_4 &\quad \sigma_4 \star pt&= q_2 D_2 - \frac{1}{2}q_2f(q_4)D_4,
\end{align*}
where $f(z)=\sum_{d\geq 1}\binom{2d}{d}z^d = \frac{1}{\sqrt{1-4z}}-1$.
\end{Thm}

\vspace{0.3cm}
We verify Theorem \ref{Thm} through a series of reductions and computations.
\subsection{Reduction to $\langle 1,D_i \rangle_{0,2,dD_4}$}\label{subsecreduce} By linearity, we only consider the monomials in $\C[\sigma_2,\sigma_4]$ for the insertion coming from the light marking. From the module structure in Theorem \ref{module}, it suffices to reduce to the case of the generators $\sigma_2$ and $\sigma_4$. Thus, it is enough to compute \[\langle T_i, T_j \mid \sigma_k \rangle_{0,2|1,\beta},\]
where $T_0:=[\F_2],~T_1:=[D_2],~T_2:=[D_4],~T_3:=[pt]$ form a basis for $H^*(\F_2)$.

In general, the forgetful map $Q_{0,m+1|k}(X,\beta) \rightarrow Q_{0,m|k}(X,\beta)$, forgetting a heavy marking, does not define a well-defined morphism in the case of the moduli space of quasimaps. On the other hand, forgetting a light marking $Q_{0,m|k+1}(X,\beta) \rightarrow Q_{0,m|k}(X,\beta)$ defines a universal curve of $Q_{0,m|k}(X,\beta)$ \cite[\S2.2]{lho}. As a consequence, the divisor equation holds for the insertion from the one light marking, i.e., 
\begin{equation}
\label{diveqn}
    \langle T_i, T_j | \sigma_k \rangle_{0,2|1,\beta}=(D_k \cdot \beta)\langle T_i, T_j \rangle_{0,2,\beta}.
\end{equation}
Therefore, proving Theorem \ref{Thm} reduces to compute 2-pointed quasimap invariants of $\F_2$ with all possible insertions for all effective curve classes $\beta \in H_2(\F_2,\Z) \simeq \Z D_2 \oplus \Z D_4$.

Write $\beta = aD_2 + bD_4$. Then, the dimension of 2-pointed quasimap moduli spaces of $\F_2$ is \[\dim Q_{0,2}(\F_2,\beta) = 1+2a.\]
Observe that one insertion from $H^*(\F_2)$ can have codimension up to 2. Since there are only 2 insertions, $a$ must be 0 or 1. The only possible pairs of insertions that we need to consider are
\begin{align*}
    \langle D_i,1 \rangle_{0,2,\beta} \quad &\text{when }\beta = dD_4,\\
    \langle D_i,pt \rangle_{0,2,\beta} \quad &\text{when }
    \beta = D_2 + dD_4,
\end{align*}
where $i=2,4$ and $d\geq 0$. Since computing $\langle D_i,pt \rangle_{0,2,\beta}$ is somewhat similar, we elaborate computation of $\langle D_i,1 \rangle_{0,2,dD_4}$ in detail.

For a component $C'$ of a quasimap, there exists a positive number $d(=\beta_{C'})$, which is the \textit{degree} of the component $C'$, such that 
\[L_1\rvert_{C'} \simeq \calO_{C'}(d),~L_2\rvert_{C'} \simeq \calO_{C'}(d),~L_3\rvert_{C'} \simeq \calO_{C'},~L_4\rvert_{C'} \simeq \calO_{C'}(-2d). \]
One key observation when $\beta=dD_4$ is that, on a component $C'$, the negative degree of the fourth line bundle forces the corresponding section of the quasimap to be the zero section along $C'$. Therefore, every component maps to $D_4$ in $\F_2$, which contains the fixed points, $p_1$ and $p_2$. In this case, it satisfies that $N_{D_4 / \F_2} \simeq O_{\PP^1}(-2)$ with $H^0(\PP^1,O_{\PP^1}(-2))=0$. This property is called \textit{rigidity} \cite{bryan}. Observe that the fact that $D_3$ and $D_4$ are disjoint divisors implies
\begin{equation}
\label{D3dD4}
    \langle D_3,1 \rangle_{0,2,dD_4}=0.
\end{equation} 
Since $D_3 = 2D_2 + D_4$ and $D_1=D_2$, there are relations
\begin{equation}
\label{relationinva}
\langle D_4,1 \rangle_{0,2,dD_4} = -2 \langle D_2,1 \rangle_{0,2,dD_4}, \quad \langle D_1,1 \rangle_{0,2,dD_4} =  \langle  D_2,1 \rangle_{0,2,dD_4}.
\end{equation}
Therefore, it is enough to compute \[\langle D_2,1 \rangle_{0,2,dD_4}.\]

\subsection{Fixed loci for $\langle D_2,1 \rangle_{0,2,dD_4}$} The $(\C^*)^4$-action on $\F_2$ induces an action on $Q_{0,2}(\F_2,\beta)$, so we can apply Atiyah--Bott localization to compute $\langle D_2,1 \rangle_{0,2,\beta}$.

Depending on existence of a base point, there are two types of components of the source curve for a quasimap in $F$. For convenience, we give the following definitions to distinguish them.
\begin{Def}
Let $F$ be a fixed locus in $Q_{0,2}(\F_2,\beta)$ under the torus action, and $(C;\underline{p};\underline{q},\underline{L},\underline{\phi})$ a quasimap in $F$. We say that a component $C' \subseteq C$ is of \textbf{base-type} if the quasimap has a base point in $C'$. Otherwise, we say $C'$ is of \textbf{map-type}. When all the components of $C$ are of map-type, we say $F$ is a fixed locus of \textbf{map-type}. Otherwise, $F$ is called a fixed locus of \textbf{base-type}.
\end{Def}
  
One can express a fixed locus in $Q_{0,2}(\F_2,\beta)$ as a decorated chain graph. The fact that the moduli spaces of quasimaps we consider have two marked points forces quasimaps to have a chain of $\PP^1$'s for their source curve. Thus, a fixed locus can correspond to a chain graph:
\begin{itemize}
    \item An edge represents one map-type component of a quasimap,
    \item A vertex between two edges represents either a node or a heavy marking.
\end{itemize}

\begin{table}[t]
\begin{tabular}{|c|c|c|c|c|} 
\hline
& $\calO^\mathcal{T}(D_1)$&$\calO^\mathcal{T}(D_2)$  & $\calO^\mathcal{T}(D_3)$ & $\calO^\mathcal{T}(D_4)$\\ \hline
$p_1=(1:0:1:0)$ & 0 & $\alpha_1-\alpha_2$ & 0 & $\alpha_3-\alpha_4-2\alpha_1$\\ \hline
$p_2=(0:1:1:0)$ & $\alpha_2-\alpha_1$ & 0& 0 & $\alpha_3-\alpha_4-2\alpha_2$\\ \hline
$p_3=(0:1:0:1)$ & $\alpha_2-\alpha_1$ & 0 & $2\alpha_2-\alpha_3+\alpha_4$& 0\\ \hline
$p_4=(1:0:0:1)$ & 0 & $\alpha_1-\alpha_2$ & $2\alpha_1-\alpha_3+\alpha_4$ & 0 \\ \hline
\end{tabular}
\caption{The weights of line bundles at each fixed points}
\label{weighttable}
\end{table}

\begin{figure}[tb] 
    \centering
\begin{tikzpicture}

\draw[very thick] (1.5,-1.5) --(1.5,1.5);
\draw[very thick] (1.5,1.5) --(-1.5,1.5);
\draw[very thick] (-1.5,1.5) --(-1.5,-1.5);
\draw[very thick] (-1.5,-1.5) --(1.5,-1.5);

\node at (-1.2,-1.2) {$p_1$};
\node at (1.2,1.2) {$p_3$};
\node at (-1.2,1.2) {$p_4$};
\node at (1.2,-1.2) {$p_2$};

\node at (0,-1.8) {$D_4$};
\node at (0,1.8) {$D_3$};
\node at (-1.8,0) {$D_2$};
\node at (1.8,0) {$D_1$};

\filldraw (-1.5,-1.5) circle (1.5pt);
\filldraw (1.5,-1.5) circle (1.5pt);
\filldraw (1.5,1.5) circle (1.5pt);
\filldraw (-1.5,1.5) circle (1.5pt);

\end{tikzpicture}
\caption{A toric diagram of $\F_2$}
\label{toricdiag}
\end{figure}
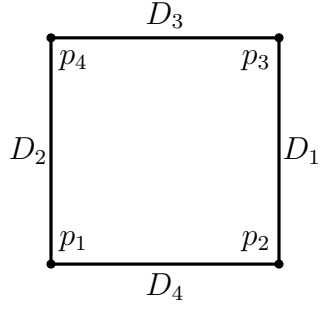

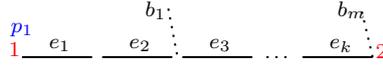
\begin{figure}[tb]
    \centering
    \tikzset{every picture/.style={line width=0.75pt}} 
\begin{tikzpicture}[x=0.75pt,y=0.75pt,yscale=-1,xscale=1]
\draw    (10,300) -- (45,300) ;
\draw    (50,300) -- (85,300) ;
\draw    (90,300) -- (125,300) ;
\draw    (150,300) -- (185,300) ;
\draw  [dash pattern={on 0.84pt off 2.51pt}]  (132,300) -- (144.6,300) ;
\draw  [dash pattern={on 0.84pt off 2.51pt}]  (83,275.6) -- (87.6,300) ;
\draw  [dash pattern={on 0.84pt off 2.51pt}]  (185,300) -- (179,274.25) ;
\draw (22,289.4) node [anchor=north west][inner sep=0.75pt]  [font=\tiny]  {$e_{1}$};
\draw (62,289.4) node [anchor=north west][inner sep=0.75pt]  [font=\tiny]  {$e_{2}$};
\draw (102,289.4) node [anchor=north west][inner sep=0.75pt]  [font=\tiny]  {$e_{3}$};
\draw (162,289.4) node [anchor=north west][inner sep=0.75pt]  [font=\tiny]  {$e_{k}$};
\draw (70,270.4) node [anchor=north west][inner sep=0.75pt]  [font=\tiny]  {$b_{1}$};
\draw (166,270.4) node [anchor=north west][inner sep=0.75pt]  [font=\tiny]  {$b_{m}$};
\draw (2,291.4) node [anchor=north west][inner sep=0.75pt]  [font=\tiny,color={rgb, 255:red, 255; green, 0; blue, 0 }  ,opacity=1 ]  {$1$};
\draw (3,281.4) node [anchor=north west][inner sep=0.75pt]  [font=\tiny,color={rgb, 255:red, 0; green, 0; blue, 255 }  ,opacity=1 ]  {$p_1$};
\draw (185,292.4) node [anchor=north west][inner sep=0.75pt]  [font=\tiny,color={rgb, 255:red, 255; green, 0; blue, 0 }  ,opacity=1 ]  {$2$};
\end{tikzpicture}
    \caption{The decorated chain graph of a general fixed locus for $\beta=dD_4$}
    \label{graphfixedlocus}
\end{figure}

We decorate our chain graph to encode the rest of the information:
\begin{itemize}
    \item A dashed half-edge at a vertex corresponds to one base-type component.
    \item A blue labelling at a vertex stands for the fixed point in the image of a quasimap.
    \item A red labelling at one end stands for a heavy marking.
    \item $e_i$ on an edge is the degree of the map-type component represented by the edge.
    \item $b_j$ on a dashed half-edge is the degree of the base-type component represented by the half-edge.
\end{itemize}

We provide informative notes on a decorated chain graph for a fixed locus:
\begin{itemize}
    \item When a dashed half-edge is attached to a vertex, the number of nodes represented by the vertex is equal to the number of edges attached to the vertex.
    \item Without loss of generality, we fix the order of the red labellings for heavy markings.
    \item We put only one blue labelling, since these will alternate along consecutive vertices by rigidity. In fact, it is redundant, but it reminds us where the marking goes.
    \item We omit the red and blue color when it is clear.
\end{itemize}

 Figure \ref{graphfixedlocus} shows the graph of a general fixed locus, including all possible types of vertices:
\begin{enumerate}
    \item a vertex at one end without any dashed half-edge,
    \item an interior vertex (i.e., not at one end) without any dashed half-edge,
    \item an interior vertex with a dashed half-edge,
    \item a vertex at one end with a dashed half-edge.
\end{enumerate}

Depending on where a vertex goes, we classify vertices in the following way: for $k \in \{1,2\}$, define
\begin{description}
    \item[$I_{k}^{m}$] the set of all interior vertices: 1) mapping to $p_k$, 2) not carrying any dashed half-edges,
    \item[$I_{k}^{b}$] the set of all interior vertices: 1) mapping to $p_k$, 2) carrying a dashed half-edge,
    \item[$I_{k}^{end}$] the set of all vertices at one end: 1) mapping to $p_k$, 2) carrying a dashed half-edge.
\end{description}
To collect vertices going to either $p_1$ or $p_2$, denote
\[I^m := I_{1}^{m} \sqcup I_{2}^{m}, \qquad I^b := I_{1}^{b} \sqcup I_{2}^{b}, \qquad I^{end} := I_{1}^{end} \sqcup I_{2}^{end},\quad I^B:=I^b \sqcup I^{end}.\]
Last, introduce the following notations to count vertices in each sets
\begin{center}
    $N_{k}^{m}:=|I_{k}^{m}|, \qquad N_{k}^{b}:=|I_{k}^{b}|, \qquad N_{k}^{end}:=|I_{k}^{end}|,~ \qquad~$\\
    $N^{m}:=|I^{m}|, \qquad N^{b}:=|I^{b}|, \qquad N^{end}:=|I^{end}|.~~~$
\end{center}

\begin{figure}[b]
    \centering
    \tikzset{every picture/.style={line width=0.75pt}} 
\begin{tikzpicture}[x=0.75pt,y=0.75pt,yscale=-1,xscale=1]

\draw    (25,220) -- (57,220) ;
\draw    (62,220) -- (94,220) ;
\draw  [dash pattern={on 0.84pt off 2.51pt}]  (59.57,219.57) -- (59.57,191.29) ;
\draw    (102.29,210.29) -- (145,210.29) ;
\draw [shift={(147,210.29)}, rotate = 180] [color={rgb, 255:red, 0; green, 0; blue, 0 }  ][line width=0.75]    (10.93,-3.29) .. controls (6.95,-1.4) and (3.31,-0.3) .. (0,0) .. controls (3.31,0.3) and (6.95,1.4) .. (10.93,3.29)   ;
\draw [shift={(100.29,210.29)}, rotate = 0] [color={rgb, 255:red, 0; green, 0; blue, 0 }  ][line width=0.75]    (10.93,-3.29) .. controls (6.95,-1.4) and (3.31,-0.3) .. (0,0) .. controls (3.31,0.3) and (6.95,1.4) .. (10.93,3.29)   ;
\draw    (156.86,219.86) -- (178.14,193.29) ;
\draw    (223,220.29) -- (207.29,194.57) ;
\draw [color={rgb, 255:red, 0; green, 11; blue, 255 }  ,draw opacity=1 ]   (215.86,199.71) -- (168.71,200) ;
\draw  [color={rgb, 255:red, 0; green, 49; blue, 255 }  ,draw opacity=1 ][fill={rgb, 255:red, 255; green, 255; blue, 255 }  ,fill opacity=1 ] (190.57,199.79) .. controls (190.57,199.12) and (191.12,198.57) .. (191.79,198.57) .. controls (192.46,198.57) and (193,199.12) .. (193,199.79) .. controls (193,200.46) and (192.46,201) .. (191.79,201) .. controls (191.12,201) and (190.57,200.46) .. (190.57,199.79) -- cycle ;
\draw    (154.29,214.29) -- (159.43,219.43) ;
\draw    (220,220.43) -- (225.29,215.71) ;

\draw (38,209.4) node [anchor=north west][inner sep=0.75pt]  [font=\scriptsize]  {$e$};
\draw (71,205.4) node [anchor=north west][inner sep=0.75pt]  [font=\scriptsize]  {$e'$};
\draw (60.57,180.69) node [anchor=north west][inner sep=0.75pt]  [font=\scriptsize]  {$b$};
\draw (51,220.83) node [anchor=north west][inner sep=0.75pt]  [font=\footnotesize]  {$v$};
\draw (62.43,220.54) node [anchor=north west][inner sep=0.75pt]  [font=\scriptsize,color={rgb, 255:red, 34; green, 0; blue, 255 }  ,opacity=1 ]  {$p_{i}$};
\draw (162,188.69) node [anchor=north west][inner sep=0.75pt]  [font=\scriptsize]  {$n$};
\draw (208.71,186.54) node [anchor=north west][inner sep=0.75pt]  [font=\scriptsize]  {$n'$};
\draw (161.14,212.97) node [anchor=north west][inner sep=0.75pt]  [font=\tiny]  {$C_{e}$};
\draw (205,213.69) node [anchor=north west][inner sep=0.75pt]  [font=\tiny]  {$C_{e'}$};
\draw (183.14,205.54) node [anchor=north west][inner sep=0.75pt]  [font=\tiny]  {$C_{b}$};
\end{tikzpicture}
    \caption{A picture of a local vertex}
    \label{localver}
\end{figure}
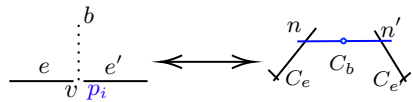

We assign to each vertex $v$ the data $(i(v),b(v),e(v),e'(v),n(v),n'(v))$, where
\begin{description}
    \item[$i(v)$] the corresponding fixed point $p_i$ to which the vertex $v$ goes,
    \item[$b(v)$] the degree of dashed half-edge attached to $v$; 0 if nothing is attached,
    \item[$e(v)$] the degree of an edge attached to $v$; this is always positive,
    \item[$e'(v)$] the degree of another line, if it exists. Otherwise, it is 0.
    \item[$n(v)$] one node that is represented by $v$
    \item[$n'(v)$] another node, if it exists.
\end{description}
We omit $v$ and write $i,b,e,e',n,n'$ when it is clear in the context. See Figure \ref{localver}.
\subsection{Virtual normal bundle}
Applying localization, we obtain a formula for the invariant $\langle D_2,1 \rangle_{0,2,dD_4}$ in terms of equivariant classes
\begin{equation}
\label{afterlocalization}
    \langle D_2,1 \rangle_{0,2,dD_4} = \sum_{F:\text{ fixed loci} } \int_F\frac{ev_1^* c_1^\mathcal{T}(\calO^\mathcal{T}(D_2))\lvert_F }{e^\mathcal{T}(N_{F}^{\text{vir}})},
\end{equation}
where $N_{F}^{\text{vir}}$ is the virtual normal bundle to $F$ in the moduli space. For a detailed explanation, we refer readers to \cite[\S9.2]{coxmirror}, \cite[\S4]{GrabPand}, or \cite[\S27.4]{mirrsymm}.

For a fixed locus $F$ in $Q_{0,2}(\F_2,\beta)$, we compute $e^\mathcal{T}(N_{F}^{\text{vir}})$. For a vector bundle $E$ with a torus action, denote its nonzero weight part by $E^{\text{mov}}$. From the tangent-obstruction sequence as in \cite[ch7 \S1.4, ch9 \S2.1]{coxmirror} or \cite[\S24.4]{mirrsymm}, the formula for the inverse of the equivariant Euler class of the virtual normal bundle to $F$ is given by
\begin{equation}
\label{virnor}
    \frac{1}{e^\mathcal{T}(N_{F}^{\text{vir}}) } = \frac{e^\mathcal{T}(\Aut(C,x_1,x_2)^{\text{mov}} )}{e^\mathcal{T}(\Defm(C,x_1,x_2)^{\text{mov}} ) ~e^\mathcal{T}(\Defm(\underline{s})^{\text{mov}} )},
\end{equation}
where $(C;\underline{x},\underline{L},\underline{s}  )$ is a quasimap in $F$.

For convenience, set notations for weights
\begin{align*}
W_1: = e^\mathcal{T}(\calO(D_2)\rvert_{p_1} )=\alpha_1 - \alpha_2,\quad V_1: = e^\mathcal{T}(\calO(D_4)\rvert_{p_1} )=\alpha_3-\alpha_4-2\alpha_1,\\
W_2: = e^\mathcal{T}(\calO(D_2)\rvert_{p_2} )=\alpha_2 - \alpha_1,\quad V_2: = e^\mathcal{T}(\calO(D_4)\rvert_{p_2} )=\alpha_3-\alpha_4-2\alpha_2.
\end{align*}
Also, we write $W:=W_1,V:=V_1$. Then, it is easy to see that
\begin{equation}
\label{vandw}
    W_1=W,~W_2=-W,~V_1=V~V_2=V+2W.
\end{equation}

We explain each factor in (\ref{virnor}). For a more detailed explanation, we refer \cite[\S27.4]{mirrsymm} to the reader.

1) Automorphisms of $(C;x_1,x_2)$: we do not have any contributions from automorphisms of a pointed source curve. The source curve of a quasimap in a fixed locus looks like a chain of $\PP^1$'s, each component with exactly two special points: a node or a heavy marking (recall that a red labelling in a decorated graph represents a heavy marking).

2) Deformation of $(C;x_1,x_2)$: a vertex $v \in I^m$ has valence 2 (dashed half-edges are not counted). Since the two map-type components attached to this vertex are not contracting, the weight from deformation of the source curve comes from the tensor product of the two tangent lines at the node in each component. Thus, the contribution to $e^\mathcal{T}(N_{F}^{\text{vir}})$ is
\begin{equation}
\label{term1}
(-1)^{N_{2}^{m}} W_1^{N^{m}} \prod_{v\in I^{m}} \bigg (\frac{1}{e} + \frac{1}{e'} \bigg) .\end{equation}

Consider that a vertex $v \in I^B$ can have one or two edges. Define 
\[\epsilon(v):= \begin{cases}0,&v\in I^{end}\\
1,&v\in I^{b}.
\end{cases}
\] For a vertex $v \in I^{end}$, there is one contracting component(dashed half-edge) and one non-contracting component(edge). On the other hand, a vertex $v \in I^{b}$ has one contracting component with two nodes to which two non-contracting components are attached. For instance, Figure \ref{localver} shows a local picture around a vertex, where the right figure is a picture of the graph in geometric side. In this case, the weight contribution from deformation of curves is obtained as follow. At the nodes $n$($n'$, resp) along $C_e$($C_{e'}$, resp), we obtain the weight $\frac{W_i}{e}$($\frac{W_i}{e'}$, resp). On the other hand, $C_b$ may be view as an element of the Losev-Manin space $\overline{M}_{2|b}$. Deformations of a contracting component gives rise to a $\psi$-class at the node, see \cite[\S25.2]{mirrsymm}. Thus, along $C_b$, we have a $\psi$-classes $\psi_n:=\psi_{1}(\overline{M}_{2|b})$ and  $\psi_{n'}:=\psi_{2}(\overline{M}_{2|b})$. The contribution to $e^\mathcal{T}(N_{F}^{\text{vir}})$ is
\begin{equation}
\label{term2}
\prod_{v \in I^{B}} \bigg(\frac{W_i}{e} - \psi_{n}  \bigg)\bigg(\frac{W_i}{e'} - \psi_{n'}  \bigg)^{\epsilon(v)} .
\end{equation}

3) Deformation of sections $\underline{s}$: There is an Euler sequence
\begin{equation}
\label{relpot}
    0\rightarrow \calO_{\mathcal{C}}^{\oplus r} \rightarrow  \oplus_\rho \calL_{\rho}\rightarrow \calF \rightarrow 0 
\end{equation} 
over a universal curve $\pi: \mathcal{C} \rightarrow Q:=Q_{0,2}(\F_2,\beta)$. In \cite[\S5.3]{toricquasi}, the relative obstruction theory over $\frak{M}_{0,2}^{\circ}$ is given by \[E_{Q/\frak{M}_{0,2}^{\circ}}^{\bullet}:=(\mathbf{R}^{\bullet}\pi_*\calF)^\vee.\]
We will use this to figure out the weight of $\text{Def}(\underline{s})$.

For the source curve $C$ of a given quasimap in a general fixed locus, one can break $C$ into its components using normalization, and the normalization exact sequence is the following:
\begin{align*}
    0\rightarrow \calO_C \rightarrow \bigoplus_{C_e:\text{map-type}}\calO_{C_e} \oplus  \bigoplus_{C_b:\text{base-type}} \calO_{C_b} \rightarrow \bigoplus_{n: \text{node}} \calO_n \rightarrow 0
\end{align*}
Tensoring $\calF$ and taking cohomology gives rise to
\begin{align*}
    0&\rightarrow H^0(C,\calF) \rightarrow \bigoplus_{C_e:\text{map-type}}H^0(C_e,\calF) \oplus  \bigoplus_{C_b:\text{base-type}} H^0(C_b,\calF) \rightarrow \bigoplus_{n: \text{node}} T_{p_n}\F_2\\
    &\rightarrow H^1(C,\calF) \rightarrow \bigoplus_{C_e:\text{map-type}}H^1(C_e,\calF) \oplus  \bigoplus_{C_b:\text{base-type}} H^1(C_b,\calF) \rightarrow 0,
\end{align*}
where $p_n$ is the fixed point in $\F_2$ where the node $n$ goes to. 

Over map-type components $C_e$, since sections define a map, it amounts to compute the weights of $H^i(C_e,(\underline{s}\lvert_{C_e})^{*} T\F_2 )$. Thus, the contribution to $1/ e^\mathcal{T}(N_{F}^{\text{vir}})$ for one map-type component of degree $e$ is 
\begin{equation}
\label{term4}
\frac{e^{2e}\prod_{j=0}^{2e-2} \big( V+ \frac{1+j}{e}W \big)
}{e
(e!)^2W^{2e}(-1)^e
},
\end{equation}
which agrees with the one from the moduli space of stable maps.

For the nodes, since $T_{p_n}\F_2 \simeq (\calO(D_i) \oplus \calO(D_4) )\lvert_{p_n}$ where $i=1,2$ depending on $p_n$, one can easily take the equivariant Euler class for their weights. The contribution to $1/ e^\mathcal{T}(N_{F}^{\text{vir}})$ is
\begin{equation}
\label{term3}
    \bigg( \prod _{v \in I^{m}} W_iV_i \bigg) \bigg( \prod _{v \in I^{end}} W_iV_i \bigg) \bigg( \prod_{v \in I^b}W_i^2 V_i^2\bigg).
\end{equation}

For a vertex $v \in I^B$, there is a corresponding base-type component $C_b$ may be regarded as an element in the Losev-Manin space $\overline{M}_{2|b}$. In this case, since $C_b$ has base points, the component is contracting to a torus fixed point \[p_1 =D_2\cap D_4 \text{ or } p_2 = D_1 \cap D_4.\]
 Thus, it is enough to consider the following from \eqref{relpot}
\[
\mathbf{R}^{\bullet}\pi_*(\mathcal{L}_2\oplus \mathcal{L}_4) \rvert_F \text{ or }\mathbf{R}^{\bullet}\pi_*(\mathcal{L}_1\oplus \mathcal{L}_4) \rvert_F,
\]
respectively. Because of the base points on $C_b$, we consider the following line bundles
\[\mathcal{O}_{C_b}(q_1 + \cdots + q_b) \text{ and } \mathcal{O}_{C_b}(-2q_1 - \cdots -2 q_b),\]
where $q_i$ are the base points on $C_b$. One can apply the divisor sequence iteratively to these:
 \begin{align*}
0 \rightarrow \mathcal{O}_{C_b}^\mathcal{T}(\sum_{k=1}^{b-1}q_k) \rightarrow \mathcal{O}_{C_b}^\mathcal{T}(&\sum_{k=1}^{b}q_k) \rightarrow  \mathcal{O}_{C_b}^\mathcal{T}(\sum_{k=1}^{b}q_k)\lvert_{q_b} \rightarrow 0\\
0 \rightarrow \mathcal{O}_{C_b}^\mathcal{T}(\sum_{k=1}^{b-2}q_k) \rightarrow \mathcal{O}_{C_b}^\mathcal{T}(&\sum_{k=1}^{b-1}q_k) \rightarrow  \mathcal{O}_{C_b}^\mathcal{T}(\sum_{k=1}^{b-1}q_k)\lvert_{q_{b-1}} \rightarrow 0\\
&\vdots\\
0 \rightarrow \mathcal{O}_{C_b}^\mathcal{T} \rightarrow \mathcal{O}_{C_b}(&q_1) \rightarrow  \mathcal{O}_{C_b}^\mathcal{T}(q_1)\lvert_{q_1} \rightarrow 0,
 \end{align*}
 and 
 \begin{align*}
    0 \rightarrow \mathcal{O}_{C_b}^\mathcal{T}(-2\sum_{k=1}^{b}q_k) \rightarrow \mathcal{O}_{C_b}^\mathcal{T}(-2\sum_{k=1}^{b-1}&q_k-q_b) \rightarrow  \mathcal{O}_{C_b}^\mathcal{T}(-2\sum_{k=1}^{b-1}q_k-q_b)\lvert_{q_b} \rightarrow 0\\
    0 \rightarrow \mathcal{O}_{C_b}^\mathcal{T}(-2\sum_{k=1}^{b-1}q_k-q_b) \rightarrow \mathcal{O}_{C_b}^\mathcal{T}&(-2\sum_{k=1}^{b-1}q_k) \rightarrow  \mathcal{O}_{C_b}^\mathcal{T}(-2\sum_{k=1}^{b-1}q_k)\lvert_{q_b} \rightarrow 0\\
    0 \rightarrow \mathcal{O}_{C_b}^\mathcal{T}(-2\sum_{k=1}^{b-1}q_k) \rightarrow \mathcal{O}_{C_b}^\mathcal{T}(-2\sum_{k=1}^{b-2}q_k&-q_{b-1}) \rightarrow  \mathcal{O}_{C_b}^\mathcal{T}(-2\sum_{k=1}^{b-1}q_k-q_{b-1})\lvert_{q_{b-1}} \rightarrow 0\\
    0 \rightarrow \mathcal{O}_{C_b}^\mathcal{T}(-2\sum_{k=1}^{b-2}q_k-q_{b-1}) \rightarrow \mathcal{O}_{C_b}^\mathcal{T}&(-2\sum_{k=1}^{b-2}q_k) \rightarrow  \mathcal{O}_{C_b}^\mathcal{T}(-2\sum_{k=1}^{b-2}q_k)\lvert_{q_{b-1}} \rightarrow 0\\
    &\vdots\\
    0 \rightarrow \mathcal{O}_{C_b}^\mathcal{T}(-2q_1) \rightarrow \mathcal{O}_{C_b}^\mathcal{T}&(-q_1) \rightarrow  \mathcal{O}_{C_b}^\mathcal{T}(-q_1)\lvert_{q_1} \rightarrow 0\\
    0 \rightarrow \mathcal{O}_{C_b}^\mathcal{T}(-q_1) \rightarrow \mathcal{O}&_{C_b}^\mathcal{T} \rightarrow  \mathcal{O}_{C_b}^\mathcal{T}\lvert_{q_1} \rightarrow 0.
 \end{align*}
Taking the long exact sequence and equivariant Euler class, one can obtain the contribution deformation of sections to $1/ e^\mathcal{T}(N_{F}^{\text{vir}})$: for $v \in I^{B}$, 
\begin{align}
\label{term5}
e^\mathcal{T}(H^0(C,&\mathcal{O}_C^\mathcal{T}(q_1 + \cdots + q_b))) ^{-1} \cdot e^\mathcal{T}(H^1(C,\mathcal{O}_C^\mathcal{T}(-2q_1 - \cdots - 2q_b)))\\
    &= \frac{1}{b!}\frac{V_i}{W_i(W_i-\hat{\psi}_1)}\prod_{j=2}^{b}\frac{(V_i-2\Delta_{j})(V_i-2\Delta_{j} + \hat{\psi}_j)}{W_i+\Delta_j}, \nonumber
\end{align}
where $D_{ij}$ is the divisor on $\overline{M}_{2|b}$ parameterizing curves with $q_i$ and $q_j$ colliding, $\Delta_i:= D_{1i} + D_{2i} + \cdots + D_{i-1,i}$, and $\hat{\psi}_j= -e^\mathcal{T}(\calO(q_j)\lvert_{q_j})$ is the $\psi$-class in $\overline{M}_{0,2|b}$ at the light point $q_j$. The factor $1/b!$ comes from permuting the base points. For more detailed explanation, we refer readers to \cite{Oprea2016}. 

Combining (\ref{term1}), (\ref{term2}), (\ref{term3}), (\ref{term4}), (\ref{term5}) together, one can write a formula of the inverse of the Euler class of the virtual normal bundle for a fixed locus $F$.
\begin{Pro}
\label{formula}
For a fixed locus $F$ of $Q_{0,2}(\F_2,dD_4)$, 
    \begin{align} \label{onestar}
    \frac{1}{e^\mathcal{T}(N_{F}^{\text{vir}}) } &=\Cont(VC)\Cont_{m}(NS)\prod_{e\in\text{Edges}}\Cont_E(e) \prod_{v \in I^{B}} \Cont_B(v),
\end{align}
where
\begin{align*}
    \Cont(VC)&:=(-1)^{N_{2}^{m}+N_{2}^{end}+2N_{2}^{b}} V_1^{N_1^{m}+N_{1}^{end}+2N_1^{b}} V_2^{N_2^{m}+N_{2}^{end}+2N_2^{b}} W_1^{N^{m}+N^{end}+2N^{b}}\\
    \Cont_{m}(NS)&:=\bigg( (-1)^{N_{2}^{m}} W_1^{N^{m}} \prod_{v\in I^{m}} \bigg (\frac{1}{e_v} + \frac{1}{e'_v} \bigg) \bigg)^{-1}\\
    \Cont_E(e)&:= \frac{1}{e}\frac{e^{2e}\prod_{j=0}^{2e-2} \big( V_1+ \frac{1+j}{e}W_1 \big)}{(e!)^2W_1^{2e}(-1)^e}\\
    \Cont_B(v)&:= \int_{\overline{M}_{0,2|b}}\frac{1}{b!}
    \frac{ \frac{V_i}{W_i^2}\prod_{j=2}^{b}\frac{(V_i-2\Delta_{j})^2}{W_i+\Delta_j} }{ \bigg(\frac{W_i}{e} - \psi_{n}  \bigg)\bigg(\frac{W_i}{e'} - \psi_{n'}  \bigg)^{\epsilon(v)}}.
\end{align*}
\end{Pro}
$\Cont(VC)$ is the contribution from vertex counting in the chain graph of $F$, $\Cont_{m}(NS)$ the contribution from node smoothing at a map-type vertex, $\Cont_E(e)$ the contribution of an edge, and $\Cont_B(v)$ the contribution of a base-type vertex. Also, note that we omitted all $\psi$-classes $\hat{\psi}_j$ from \eqref{psiformulawithlight}.

\subsection{Simplifying the formula} We simplify the formula \eqref{onestar} to obtain a rational number through a series of observations. One can observe that the denominator in \eqref{onestar} would be of the form $W^N$. This allows us to consider some fixed loci that contribute to the final answer. As a result, we do not need to figure out the complete expansion of $\Cont_B(v)$ which requires a somewhat complicated combinatorics problem.

\begin{Lem}
\label{denomWk}
    The denominator in \eqref{onestar} is of the form $W^N$ for some $N$.
\end{Lem}
\begin{proof}
    It is clear by applying geometric series to the term $\big(\frac{W_i}{e} - \psi_{n}  \big)\big(\frac{W_i}{e'} - \psi_{n'}  \big)$. Recall that $W_i  = \pm W$.
\end{proof}
\begin{figure}[b]
\tikzset{every picture/.style={line width=0.75pt}} 
\begin{tikzpicture}[x=0.75pt,y=0.75pt,yscale=-1,xscale=1]
\draw    (60,68) -- (97,68) ;
\draw    (128,68) -- (152,68) ;
\draw    (158,68) -- (182,68) ;
\draw    (60,109) -- (97,109) ;
\draw    (128,109) -- (152,109) ;
\draw    (158,109) -- (182,109) ;
\draw  [dash pattern={on 0.84pt off 2.51pt}]  (85,90) -- (100,109) ;
\draw  [dash pattern={on 0.84pt off 2.51pt}]  (140,90) -- (155,109) ;
\draw (74,56) node [anchor=north west][inner sep=0.75pt]  [font=\tiny]  {$d$};
\draw (135,59) node [anchor=north west][inner sep=0.75pt]  [font=\tiny]  {$e$};
\draw (157,56) node [anchor=north west][inner sep=0.75pt]  [font=\tiny]  {$d-e$};
\draw (135,100) node [anchor=north west][inner sep=0.75pt]  [font=\tiny]  {$e$};
\draw (164,97) node [anchor=north west][inner sep=0.75pt]  [font=\tiny]  {$e'$};
\draw (74,100) node [anchor=north west][inner sep=0.75pt]  [font=\tiny]  {$e$};
\draw (88,84) node [anchor=north west][inner sep=0.75pt]  [font=\tiny]  {$b$};
\draw (145,85) node [anchor=north west][inner sep=0.75pt]  [font=\tiny]  {$b$};

\draw (30,66) node [anchor=west][inner sep=0.75pt]  [font=\tiny]  {($F_d$)};
\draw (185,66) node [anchor=west][inner sep=0.75pt]  [font=\tiny]  {($F_{e,d-e}$)};
\draw (30,105) node [anchor=west][inner sep=0.75pt]  [font=\tiny]  {($F_d^b$)};
\draw (185,105) node [anchor=west][inner sep=0.75pt]  [font=\tiny]  {($F_{e,e'}^{b}$)};
\end{tikzpicture}
\caption{Necessary fixed loci for $\beta=dD_4$}
\label{reducedloci}
\end{figure}

For a fixed degree $\beta=dD_4$, there are some fixed loci enough to consider to compute $\langle D_2,1 \rangle_{0,2,dD_4}$. We call such a fixed locus a \textbf{necessary fixed locus}. The following corollary classifies necessary fixed loci.
\begin{Cor}
\label{correducedloci}
The graphs of necessary fixed loci for $\beta=dD_4$ are those in Figure \ref{reducedloci}, where $b+e+e'=d$. Write it as $F_{e,e'}^{b}$, and we omit indices if they are 0.
\end{Cor}
\begin{proof}
For each fixed locus $F$, the insertion at the first/second marked point gives rise to the numerator \[ev_1^*c_1^\mathcal{T}(\mathcal{O}^\mathcal{T}(D_2))\lvert_F = W,\]
since the image of a quasimap in $F$ lies on $D_4$. On the other hand, Lemma \ref{denomWk} tells us that for each fixed locus $F$, the denominator in \eqref{onestar} is $W^N$ for some $N$. Since the invariant $\langle D_2,1 \rangle_{0,2,dD_4}$ generates a rational number, the factor
\[V_1^{N_1^{m}+N_{1}^{end}+2N_1^{b}}\]
in \eqref{onestar} must vanish to contribute to the final answer. Hence, necessary fixed loci have $N_1^{m}=0,~N_{1}^{end}=0,~N_1^{b}=0$. The four graphs in Figure \ref{reducedloci} are all the graphs of the necessary fixed loci satisfying this requirement.
\end{proof}
\begin{Rem}
This corollary tells us that a vertex with a dashed half-edge must contracts to $p_2$. Thus, we will set up $i=2$ in $V_i$ and $W_i$ for such a vertex.
\end{Rem}

To reduce \eqref{onestar}, it is necessary to manipulate further the denominator and the numerator of $\Cont_B(v)$ in \eqref{onestar} for a vertex $v\in I^B$. Recall that  $\Cont_B(v)$ is given as an integration over the Losev-Manin space $\overline{M}_{0,2|b}$ whose $\dim$ is $b-1$.

i) (Denominator) The inverse of \eqref{term2}, when $\epsilon(v)=1$, can be written as:
\begin{align*}
&\frac{ee'}{W_i^2}\int_{\overline{M}_{0,2|b}}\sum_{s=0}^{b-1}\bigg(\frac{e\psi_n}{W_i}\bigg)^s\sum_{t=0}^{b-1}\bigg(\frac{e'\psi_{n'}}{W_i}\bigg)^t\\
&=\frac{ee'}{W_2^2}\int_{\overline{M}_{0,2|b}}\sum_{k=0}^{b-1}\frac{1}{W_2^k} \sum_{m=0}^{k}(e\psi_n)^{k-m}(e'\psi_{n'})^m\\
&= \frac{ee'}{W^2}\int_{\overline{M}_{0,2|b}}\sum_{k=0}^{b-1}\frac{(-1)^k}{W^k} \sum_{m=0}^{k}(e\psi_n)^{k-m}(e'\psi_{n'})^m.
\end{align*}
ii) (Numerator) Recall the a vertex $v\in I^b$ must go to $p_2$, so take $i=2$. We can write
\begin{equation}
\label{Vzero}
    \frac{V_i}{W_i^2}\prod_{j=2}^{b}\frac{(V_i-2\Delta_{j})^2}{W_i+\Delta_j} = \frac{V_2}{W_2^2}\prod_{j=2}^{b}\frac{(V_2-2\Delta_{j})^2}{W_2+\Delta_j} = \frac{V+2W}{W^2}\prod_{j=2}^{b}\frac{(V+2W-2\Delta_{j})^2}{-W+\Delta_j}.
\end{equation}
Since a term in the numerator containing $V$ does not contribute to the final answer, we take $V=0$ in \eqref{Vzero}. Then, we obtain
\begin{equation*}
    \frac{(-1)^{b-1}2^{2b-2}2W}{W^2}\prod_{j=2}^{b}\frac{(W-\Delta_{j})^2}{W-\Delta_j} =\frac{(-1)^{b-1}2^{2b-1}}{W}\prod_{j=2}^{b}(W-\Delta_{j}) .
\end{equation*}
Thus, assuming $\epsilon(v)=1$, i.e., $v\in I^b$, 
\begin{align}
\label{simplifythis}
 \Cont_B(v)\lvert_{V=0}&=\frac{(-1)^{b-1} ee' 2^{2b-1}W^{b-4}}{b!}\cdot\\
 &\int_{\overline{M}_{0,2|b}}
    \prod_{j=2}^{b}(1-\frac{\Delta_{j}}{W}) \sum_{k=0}^{b-1}\frac{1}{(-W)^k} \sum_{m=0}^{k}(e\psi_n)^{k-m}(e'\psi_{n'})^m \nonumber
\end{align}

\subsection{Calculus on $\overline{M}_{0,2|b}$}\label{comb1} To simplify \eqref{simplifythis}, we give some facts about calculus on $\overline{M}_{0,2|b}$.

For the intersection of $\psi$-classes at heavy points, it is known from \cite[Example 4.5]{MOP} that
\begin{equation}
\label{psiformula}
\int_{\overline{M}_{0,2|b}}\psi_{n}^{b-1-m}\psi_{n'}^m = \binom{b-1}{m}.
\end{equation}
Also, observe that the locus corresponding to $D_{\{i_1,\ldots,i_l \}}$ is naturally isomorphic to as the Losev-Manin space $\overline{M}_{0,2|b-l+1}$. It follows that
\begin{equation}
\label{reducelights}
    \int_{\overline{M}_{0,2|b}} D_{\{i_1,\ldots,i_{j+1} \}}
    \psi_{n}^{b-1-j-m}\psi_{n'}^{m} =  \int_{\overline{M}_{0,2|b-j}}
    \psi_{n}^{b-1-j-m}\psi_{n'}^{m}
\end{equation}

On the other hand, the intersection of $\psi$-classes at both heavy and light points is the following from \cite[\S4.6]{MOP}
\begin{equation}
\label{psiformulawithlight}
\int_{\overline{M}_{0,2|b}}\prod_{i=1}^{2}\psi_{i}^{n_i}\prod_{j=1}^{b}\hat{\psi}_{j}^{m_j}=0
\end{equation}
if $m_j\neq 0$ for some $j$.

There are classes coming from collisions of light markings.
\begin{Def}
Let $b$ be a positive integer and $1\leq l\leq b$. For a subset $I=\{i_j\}_{j=1}^{l} \subseteq [b]$ of $l$ distinct elements, define \[D_{\{i_1,\ldots,i_l \}} \] to be the cycle class in $H^*(\overline{M}_{0,2|b})$ of the closure of the locus where all $i_j$th light markings are colliding. Similarly, for $I' \subseteq [b]$ with $I \cap I' =\varnothing$,  denote \[D_I D_{I'} \in H^*(\overline{M}_{0,2|b}) \] by the cycle class for the closure of the locus where light markings in $I$ are colliding with the ones in $I$, but not any markings in $I'$, and vice versa.
\end{Def}
When $l=1$, $D_{\{i_1\}}=[\overline{M}_{0,2|b}]=:1$. The codimension of $D_{\{i_1,\ldots,i_l\}}$ is $l-1$. We provide a simple lemma for the product of two such cycle classes in $H^*(\overline{M}_{0,2|b})$.
\begin{Lem}
\label{productD}
    For subsets $I,I' \subseteq [b]$,
    \[D_I \cdot D_{I'} = \begin{cases}
    D_{I}D_{I'},  &\text{if }|I \cap I' |=0\\
    D_{I \cup I' },  &\text{if }|I \cap I' |=1 \\
    D_{I \cup I' }(-\hat{\psi}_{I\cap I'})^{|I\cap I'| - 1}, &\text{if }|I \cap I' |>1,
\end{cases}\]
where $\hat{\psi}_{J}:=\hat{\psi}_{j} \lvert_{D_J}$ for every $j \in J$ (note that $\hat{\psi}_{j} \lvert_{D_J}=\hat{\psi}_{j'} \lvert_{D_J}$ for every $j,j'\in J$).
\end{Lem}
\begin{proof}
    When $|I \cap I' |=1$, $D_I \cdot D_{I'}$ is the locus where all light marking in $I\cap I'$ are colliding without any repetition. Thus, $D_I \cdot D_{I'} = D_{I \cup I' }$.

    Write $D_{12}:= D_{\{1,2\}}$ and $-\hat{\psi}_{12}:= -\hat{\psi}_{\{1,2\}}$. We claim \begin{equation}
    \label{selfintD12}
    D_{12}^2=(-\hat{\psi}_{12})D_{12}.    
    \end{equation}
    Recall the fact that $D_{12} \simeq \overline{M}_{0,2|1} \xhookrightarrow{i}  \overline{M}_{0,2|2}$. Consider the exact sequence
    \[T_{\overline{M}_{0,2|1}} \hookrightarrow i^* T_{\overline{M}_{0,2|2}} \twoheadrightarrow N_i. \]
    Then, since $\overline{M}_{0,2|1} \simeq \text{pt}$,
    \[D_{12}^2=D_{12} \cdot e(N_i)=  -D_{12} \cdot e( T_{D_{12}}^{*}) = -D_{12} \cdot \hat{\psi}_{12}. \]

    Observe that for $I=\{1,2,\ldots,b\}$, using \eqref{selfintD12} gives
    \begin{align}
    \label{selfintDI}
        D_I^2 = (D_{12} \cdot D_{23} \cdot \cdots \cdot D_{b-1,b})^2 = (-\hat{\psi}_{12\cdots b})^{b-1}(D_{12} \cdot D_{23} \cdot \cdots \cdot D_{b-1,b})=(-\hat{\psi}_{I})^{b-1} D_I.
    \end{align}
    Thus, in general, for subsets $I,I' \subseteq [b]$ with $|I \cap I' |>1$ and $i \in I \cap I'$, applying \eqref{selfintDI} implies 
    \begin{align*}
        D_I \cdot D_{I'} &= D_{(I\backslash (I \cap I'))\cup \{i\}} \cdot D_{(I'\backslash (I \cap I'))\cup \{i\}} \cdot D_{I\cap I'}^2\\
        &= D_{(I\backslash (I \cap I'))\cup \{i\}} \cdot D_{(I'\backslash (I \cap I'))\cup \{i\}} \cdot D_{I\cap I'} (-\hat{\psi}_{I\cap I'})^{|I\cap I'| - 1}\\
        &= D_{I \cup I'} (-\hat{\psi}_{I\cap I'})^{|I\cap I'| - 1}
    \end{align*}

    When $|I\cap I'|=0$, we claim that $D_I$ and $D_{I'}$ are not equivalent. Without loss of generality, assume that $|I|=|I'|$. Otherwise, codimensions of $D_I$ and $D_{I'}$ are different. Let $b:=|I|+|I'|$. Without loss of generality, we may assume $I=\{1,\ldots, k\},~ I'=\{k+1,\ldots, b\}$. Suppose that $D_I$ and $D_{I'}$ are equivalent. Then, observe
    \begin{align*}
        D_I \cdot D_{I'} \cdot D_{k,k+1} &= D_I^2 \cdot D_{k,k+1} \\
        D_{12\cdots b} &= (-\hat{\psi}_{I})^{k-1}D_{12\cdots k+1}.
    \end{align*}
    Taking the integration over $\overline{M}_{0,2|b}$ leads us to a contradiction with \eqref{psiformulawithlight}
    \[\int_{\overline{M}_{0,2|b}}D_{12\cdots b} = \int_{\overline{M}_{0,2|1}}1 =1 \neq 0=\int_{\overline{M}_{0,2|b}}(-\hat{\psi}_{I})^{k-1}D_{12\cdots k+1} = \int_{\overline{M}_{0,2|b-k}}(-\hat{\psi}_{I})^{k-1} .\]
\end{proof}

\subsection{Unordered set partitions}\label{comb2} We also introduce basic combinatorics to reduce \eqref{simplifythis}. Let $b$ be a positive integer.
Denote  by $P(b)$ the set of all unordered set partitions of $[b]$. We always write an element $P \in P(b)$ as $P=(P_1,\ldots, P_l)$ such that $p_1:=|P_1|\geq p_2:=|P_2|\geq \cdots \geq p_l:=|P_l|$. Also, for $\lambda \vdash b $, denote $P(b,\lambda)$ by the subset of $P(b)$ satisfying $p_i = \lambda_i$ for all $i=1,2,\ldots, l=l(\lambda)$, where $l(\lambda)$ is the length of $\lambda$. There is a decomposition of $P(b)$
\[P(b) = \bigsqcup_{\lambda \vdash b} P(b,\lambda).\]
For a partition $\lambda \vdash b$, the size $|P(b,\lambda)|$ is known
\[|P(b,\lambda)| = \frac{b!}{\prod_{N=1}^{\infty}(N!)^{k_N}k_N!   } ,\]
where $(k_N)_{N=1}^{\infty}$ is the multiplicity sequence for $\lambda$ by letting $k_N$ be the number of $N$'s in $\lambda$.

We introduce notations for the class $D_I$ related to the above unordered set partitions.
If $P=(P_1,\ldots, P_l)\in P(b)$, write
\[D_P:= D_{P_1}D_{P_2}\cdots D_{P_l}.\]
The codimension of $D_P$ is $b-l$. For $\lambda \vdash b$,
\begin{align*}
    D_\lambda&:= \sum_{ P \in P(b,\lambda) } D_P.
\end{align*}
The following is a lemma for simplifying \eqref{simplifythis}.
\begin{Lem} The following identities hold in $H^*(\overline{M}_{0,2|b})$:
\label{expansion}
\begin{align*}
    \prod_{j=2}^{b}(1+\Delta_{j}) = \sum_{\lambda \vdash b} \bigg(\prod_{q=1}^{l(\lambda)}(\lambda_q-1)! \bigg) D_\lambda = \sum_{l=1 }^{b} \sum_{\lambda \vdash b: l(\lambda)=l}\bigg(\prod_{q=1}^{l(\lambda)}(\lambda_q-1)! \bigg) D_\lambda.
\end{align*}
\end{Lem}
\begin{proof} We formally expand the left-hand side, apply Lemma \ref{productD}, and then we count the number of $D_I$ for a subset $I \subseteq [b]$.

For $\lambda \vdash b$ and $P \in P(b,\lambda)$, there is a way to produce $D_P$ by choosing $1$ or $D_{ij}$ in each factor $(1+\Delta_j)$ and multiplying. Thus, one can write
    \[\prod_{j=2}^{b}(1+\Delta_{j}) =\sum_{\lambda \vdash b} \sum_{P \in P(b,\lambda)} c_P D_P .\]
    Since $P_k$ are disjoint, the number of ways to form each $D_{P_k}$ is independent of $P_k$'s, say $c_{P_k}$. It is possible to write \[c_P D_P = \prod_{k=1}^{l(\lambda)}c_{P_k}D_{P_k}.\] 
    We count $c_{P_k}$. Without loss of generality, we may write $P_k = \{1,2,\ldots, \lambda_k\}$. First, we choose $D_{i\lambda_k}$ in the factor $1+\Delta_{\lambda_k}$, which is $|P_k \backslash \{\lambda_k\}|$ amount of choices. Next, we have $|P_k \backslash \{\lambda_k,\lambda_k -1\}|$ amount of choices from $1+\Delta_{\lambda_k -1}$. This way amounts to 
    \[
    c_{P_k} = (\lambda_k-1)!.
    \]
    Therefore, it proves the first equality.

    The second equality follows from expanding the first expansion with respect to codimensions. 
\end{proof}

Applying the second identity in Lemma \ref{expansion} to the factor $\prod_{j=2}^{b}(1-\frac{\Delta_{j}}{W})$,  one can write \eqref{simplifythis} as follows: 
\begin{align}\label{happy}
&\frac{(-1)^{b-1} ee' 2^{2b-1}W^{b-4}}{b!} \cdot \\ &\sum_{l=1 }^{b} \sum_{\lambda \vdash b: l(\lambda)=l} \sum_{k=0}^{b-1} \sum_{m=0}^{k}  
    \bigg(\prod_{q=1}^{l(\lambda)}(\lambda_q-1)! \bigg) \frac{(-1)^{b-l+k}e^{k-m}e'^m}{W^{b-l+k}}  \int_{\overline{M}_{0,2|b}} D_\lambda \psi_n^{k-m} \psi_{n'}^{m} \nonumber
\end{align}
To obtain nonzero values for the integration, $\dim \overline{M}_{0,2|b} =b-1$ must be the same as
\[\text{codim}\:(D_\lambda \psi_n^{k-m} \psi_{n'}^{m}) = (b-l)+ (k-m) + m = b-l+ k.\]
Set $k:=l-1$. Then, we use \eqref{psiformula} and \eqref{reducelights} to reduce \eqref{happy}, and obtain
\begin{equation}
\label{withoutfactor}
\Cont_B(v) = \frac{ ee' 2^{2b-1}}{b!W^{3}} \cdot \sum_{l=1 }^{b}  (e+e')^{l-1} \sum_{\lambda \vdash b: l(\lambda)=l} \bigg(\prod_{q=1}^{l(\lambda)}(\lambda_q-1)! \bigg) |P(b,\lambda)| 
\end{equation}
\subsection{Symmetric functions theory}\label{comb3} To reduce \eqref{withoutfactor}, we need basic symmetric function theory. 
We follow notations from \cite{stanley}. Denote $p(b,\lambda)$ to be the subset of permutations in $S_b$ whose cycle type is $\lambda$. Also, let \[z_\lambda:= \prod_{N=1}^{\infty} (N)^{k_{N}}k_{N}!,\]
where $(k_N)_{N=1}^{\infty}$ is the multiplicity sequence for $\lambda$ as before.
Then, one can easily see
\begin{equation}
\label{orderunorder}
     \bigg(\prod_{q=1}^{l(\lambda)}(\lambda_q-1)! \bigg) |P(b,\lambda)| = |p(b,\lambda)| = b! z_{\lambda}^{-1}.
\end{equation}
For positive numbers $b$ and $e$, introduce the multiset coefficient
\[\multiset{b}{e}:= \binom{b+e-1}{e},\]
which counts the number of monomials of degree $e$ in $b$ variables. One can have the following expression for a multiset coefficient.
\begin{Lem}
\label{multicoef}
    For positive numbers $b$ and $e$, the multiset coefficient $\multiset{b}{e}$ can be written as follows:
    \begin{equation}
        \multiset{b}{e} = \sum_{l=1}^{b}\frac{e^{l-1}}{(b-1)!}\sum_{\lambda \vdash b: l(\lambda)= l } |p(b,\lambda)|
    \end{equation}
\end{Lem}
Before starting a proof, recall that homogeneous and power sum symmetric functions, denoted by $h_b$ and $p_\lambda$, are the following:
\begin{align*}
    h_b&:= \sum_{1\leq i_1\leq \cdots \leq i_b\leq b} x_{i_1}\cdots x_{i_b}\\
    p_k&:=\sum_{i=1}^{b}x_{i}^{k},\quad 
    p_{\lambda}:=p_{\lambda_1}\cdots p_{\lambda_l}.
\end{align*}
\begin{proof}
    Proposition 7.7.6 in \cite{stanley} gives us
    \[h_b = \sum_{\lambda \vdash b}z_{\lambda}^{-1}p_\lambda.\]
    Evaluating $(1,\ldots,1)$ of length $e$ to both $h_b$ and $p_\lambda$ above and using the second equality in \eqref{orderunorder} give us
    \begin{align*}
    h_b(1,\ldots,1) &= \sum_{\lambda \vdash b}z_{\lambda}^{-1}p_\lambda(1,\ldots,1) = \sum_{l=1}^{b}\sum_{\lambda \vdash b:l(\lambda)=l}\frac{|p(b,\lambda)|}{b!}e^l\\
    \multiset{e}{b} &=\frac{e}{b} \sum_{l=1}^{b}\frac{e^{l-1}}{(b-1)!}\sum_{\lambda \vdash b: l(\lambda)= l } |p(b,\lambda)|\\ \multiset{b}{e} &= \sum_{l=1}^{b}\frac{e^{l-1}}{(b-1)!}\sum_{\lambda \vdash b: l(\lambda)= l } |p(b,\lambda)|.
    \end{align*}
\end{proof}
We can obtain a concise simplification for \eqref{withoutfactor}.
 \begin{Cor}
 \label{factorvinIB}
For a necessary fixed locus $F$, assuming $\epsilon(v)=1$, and a vertex $v \in I^B$,
\begin{align*}
    \Cont_B(v) \lvert_{V=0}= \frac{ ee' 2^{2b-1}}{bW^{3}} \multiset{b}{e+e'}.
\end{align*}
\end{Cor}
\begin{proof}
    Use the first equality in \eqref{orderunorder} to \eqref{withoutfactor}, then apply Lemma \ref{multicoef}.
\end{proof}

\subsection{Computation of the invariants $\langle D_2,1 \rangle_{0,2,dD_4}$} We complete our computation of the invariants $\langle D_2,1 \rangle_{0,2,dD_4}$. The following is a contribution to a term in \eqref{afterlocalization} for each necessary fixed locus.
\begin{Pro}
\label{reformulation}
Let $F$ be a necessary fixed locus and $b+e+e'=d$. Then, the contribution of the term corresponding to $F$ in \eqref{afterlocalization} to the invariant $\langle D_2,1 \rangle_{0,2,dD_4}$ is the following:
\begin{equation} \label{jaenumberyourequantions}
   \frac{ev_1^* c_1^\mathcal{T} ( \mathcal{O}^{\mathcal{T}}(D_2))\lvert_{F}}{e^\mathcal{T}(N_{F}^{\text{vir}}) }\bigg\lvert_{V=0}
   = \begin{dcases*}
        \frac{(-1)^d}{2d}\binom{2d}{d}, &$F=F_d$;\\
       \frac{(-1)^d}{2d}\binom{2e}{e}\binom{2e'}{e'}, &$F=F_{e,e'}$;\\
       \frac{(-1)^{d-b} 2^{2b}}{2b}\binom{d-1}{b-1}\binom{2e}{e}, &$F=F_{e}^{b}$;\\
       \frac{(-1)^{d-b} 2^{2b}}{2b}\binom{d-1}{b-1} \binom{2e}{e}\binom{2e'}{e'}, &$F=F_{e,e'}^{b}$.
   \end{dcases*}
\end{equation}
\end{Pro}
\begin{proof}
As in the proof of Corollary \ref{correducedloci}, the numerator is $\displaystyle ev_1^* c_1^\mathcal{T} ( \mathcal{O}^{\mathcal{T}}(D_2))\lvert_F =W$. Once we derive the formulas in \eqref{jaenumberyourequantions} when $F=F_{e,e'},~F_{e,e'}^{b}$, the rest cases can be done in a similar fashion. Assume $F=F_{e,e'}$. Using Proposition \ref{formula} with Corollary \ref{factorvinIB} to $F$, one can derive the following
\begin{align*}
    \frac{W}{e^\mathcal{T}(N_{F}^{\text{vir}}) }\bigg\lvert_{V=0} &=W\frac{e^{2e}\prod_{j=0}^{2e-2} \big( \frac{1+j}{e}W \big)}{e(e!)^2W^{2e}(-1)^e} \cdot \frac{e'^{2e'}\prod_{j=0}^{2e'-2} \big( \frac{1+j}{e'}W \big)}{e'(e'!)^2W^{2e'}(-1)^{e'}} \cdot \frac{ (-1) (2W) W }{(-1) W  \bigg (\frac{1}{e} + \frac{1}{e'} \bigg)} \\
    & = W\cdot \frac{(-1)^{e}}{2eW}\binom{2e}{e} \frac{(-1)^{e'}}{2e'W}\binom{2e'}{e'}\cdot \frac{2W^2 ee'}{W(e+e')} = \frac{(-1)^{d} }{2d} \binom{2e}{e} \binom{2e'}{e'}.
\end{align*}
Next, let $F=F_{e,e'}^{b}$. Similarly, we obtain
\begin{align*}
    \frac{W}{e^\mathcal{T}(N_{F}^{\text{vir}}) }\bigg\lvert_{V=0} &=W\frac{e^{2e}\prod_{j=0}^{2e-2} \big( \frac{1+j}{e}W \big)}{e(e!)^2W^{2e}(-1)^e} \cdot \frac{e'^{2e'}\prod_{j=0}^{2e'-2} \big( \frac{1+j}{e'}W \big)}{e'(e'!)^2W^{2e'}(-1)^{e'}} \\
    &\cdot \frac{(-1)^{2} (2W)^{2} W^{2 }}{1} \cdot  \frac{ee'2^{2b-1}}{bW^{3}} \multiset{b}{e+e'}  \\
    & =W \frac{(-1)^{e+e'}}{4ee'W^2}\binom{2e}{e}\binom{2e'}{e'} \cdot (-1)^2 4W^4 \cdot \frac{ee' 2^{2b-1}}{bW^3}\multiset{b}{e+e'} \\
    &=\frac{(-1)^{d-b} 2^{2b}}{2b}\binom{d-1}{b-1} \binom{2e}{e}\binom{2e'}{e'}.
\end{align*}
\end{proof}
Finally, we sum up all the cases in Proposition \ref{reformulation} to end our computation for the 2-pointed quasimap invariant $\langle D_2,1 \rangle_{0,2,dD_4}$.
\begin{Pro}
\label{computedD4} All 2-pointed degree $dD_4$ quasimap invariants of $\F_2$ are given by
\begin{align*}
\langle D_1,1 \rangle_{0,2,dD_4} &= \langle D_2,1 \rangle_{0,2,dD_4} = -\frac{1}{2d}\binom{2d}{d},\\
\langle D_3,1 \rangle_{0,2,dD_4}&=0,\quad \langle D_4,1 \rangle_{0,2,dD_4} = \frac{1}{d}\binom{2d}{d}.
\end{align*}
\end{Pro}
\begin{proof}
By the Atiyah--Bott localization theorem the invariant $\langle D_2,1 \rangle_{0,2,dD_4}$ is obtained by summing the contributions of each necessary fixed locus with setting $V=0$:
\begin{align*}
    \langle D_2,1 \rangle_{0,2,dD_4} &= \frac{ev_{1}^{*}c_{1}^\mathcal{T}(\mathcal{O}^\mathcal{T}(D_2))\lvert_{F_d}}{e^\mathcal{T}(N_{F_d}^{\text{vir}}) }\bigg\lvert_{V=0}
    + \sum_{e=1}^{d-1} \frac{ev_{1}^{*}c_{1}^\mathcal{T}(\mathcal{O}^\mathcal{T}(D_2))\lvert_{F_{e,d-e}}}{e^\mathcal{T}(N_{F_{e,d-e}}^{\text{vir}}) }\bigg\lvert_{V=0} \\
    +&\sum_{b=1}^{d-1} \frac{ev_{1}^{*}c_{1}^\mathcal{T}(\mathcal{O}^\mathcal{T}(D_2))\lvert_{F_{d-b}^{b}}}{e^\mathcal{T}(N_{F_{d-b}^{b}}^{\text{vir}}) }\bigg\lvert_{V=0}
    +\sum_{b=1}^{d-2}\sum_{e=1}^{d-b-1}\frac{ev_{1}^{*}c_{1}^\mathcal{T}(\mathcal{O}^\mathcal{T}(D_2))\lvert_{F_{e,d-b-e}^{b}}}{e^\mathcal{T}(N_{F_{e,d-b-e}^{b}}^{\text{vir}}) }\bigg\lvert_{V=0}.
\end{align*}
Apply Proposition \ref{reformulation}, and observe that the first term is in fact a summand of the first summation as the case $e=0$. Then,
\begin{align*}
    & \sum_{e=0}^{d-1} \frac{(-1)^d}{2d}\binom{2e}{e}\binom{2(d-e)}{d-e}
    +\sum_{b=1}^{d-1} \frac{(-1)^{d-b} 2^{2b}}{2b}\binom{d-1}{b-1}\binom{2(d-b)}{d-b} \\
    +&\sum_{b=1}^{d-2} \sum_{e=1}^{d-b-1}\frac{(-1)^{d-b} 2^{2b}}{2b}\binom{d-1}{b-1} \binom{2e}{e}\binom{2(d-b-e)}{d-b-e}.
\end{align*}
From the second summation, we leave the case $b=d-1$ which is $-4^{d-1}$, and the rest can be inside of the last summation as the case $e=0$ since $\frac{1}{b}\binom{d-1}{b-1} = \frac{1}{d}\binom{d}{b}$. Then, 
\begin{align*}
    \sum_{e=0}^{d-1} \frac{(-1)^d}{2d}\binom{2e}{e}\binom{2(d-e)}{d-e} - 4^{d-1} +\sum_{b=1}^{d-2} \frac{(-1)^{d-b} 2^{2b}}{2d}\binom{d}{b}\sum_{e=0}^{d-b-1}\binom{2e}{e}\binom{2(d-b-e)}{d-b-e}.
\end{align*}
The first and second terms can collapse into the last summation as the case $b=0$ and $b=d-1$, respectively. Thus,
\begin{equation}
\label{eqnabove2}
\sum_{b=0}^{d-1} \frac{(-1)^{d-b} 2^{2b}}{2d}\binom{d}{b}\sum_{e=0}^{d-b-1}\binom{2e}{e}\binom{2(d-b-e)}{d-b-e}.
\end{equation}
Using the formula $\displaystyle \sum_{k=0}^{n}\binom{2k}{k}\binom{2(n-k)}{n-k}=4^n$, one can write the equation \eqref{eqnabove2} as
\begin{equation}
\label{eqnabove3}
    \sum_{b=0}^{d-1} \frac{(-1)^{d-b} 2^{2b}}{2d}\binom{d}{b}\bigg( 4^{d-b} -\binom{2(d-b)}{d-b} \bigg).
\end{equation}
The term $\sum_{b=0}^{d-1} \frac{(-1)^{d-b} 2^{2b}}{2d}\binom{d}{b} 4^{d-b}$ is in fact the case $b=d$ of $\sum_{b=0}^{d-1} \frac{(-1)^{d-b-1} 2^{2b}}{2d}\binom{d}{b}\binom{2(d-b)}{d-b}$. Thus, the equation \eqref{eqnabove3} becomes
\begin{align*}
\sum_{b=0}^{d} \frac{(-1)^{d-b-1} 4^{b}}{2d}\binom{d}{b}\binom{2(d-b)}{d-b}.
\end{align*}
Apply the formula $\binom{2n}{n} = (-4)^n\binom{-1/2}{n}$ to have
\begin{align*}
    -\frac{4^{d}}{2d}\sum_{b=0}^{d}\binom{d}{b}\binom{-1/2}{d-b}.
\end{align*}

For the last step to obtain $-\frac{1}{2d}\binom{2d}{d}$, it is enough to show that the coefficient of $x^d$ from the following is $4^{-d}\binom{2d}{d}$
\begin{align*}
    (1+x)^{d-1/2}&=\sum_{l=0}^{d}x^l\sum_{b=0}^{l}\binom{d}{b}\binom{-1/2}{d-b} + (\text{higher order terms}).
\end{align*}
To achieve this, write \[(1+x)^{d+1-1/2} = \cdots + a_{d+1}x^{d+1} +\cdots.\] Observe that the induction hypothesis for $d$ allows us to write the derivative of $(1+x)^{d+1-1/2}$ as
\[(d+1/2)(1+x)^{d-1/2} = (d+1/2) (\cdots + 4^{-d}\binom{2d}{d}x^d + \cdots).\] Therefore, $a_{d+1}$ is equal to \[ \frac{1}{d+1} (d+1/2)4^{-d}\binom{2d}{d} =  4^{-(d+1)}\binom{2(d+1)}{d+1}.\] 
We computed $\langle D_3,1 \rangle_{0,2,dD_4}=0$ from \eqref{D3dD4}. From the relation \eqref{relationinva}, we can derive $\langle D_4,1 \rangle_{0,2,dD_4}=\frac{1}{d}\binom{2d}{d}$.
\end{proof}
Therefore, we successfully computed $\langle D_i,1 \rangle_{0,2,dD_4}$ for all $i=1,2,3,4$.

\subsection{Computation of the invariants $\langle D_i,pt \rangle_{0,2,\beta}$ for $\beta=D_2 + dD_4$} To have a full description of the quantum module structure for $\F_2$, it still remains to compute 2-pointed quasimap invariants for degree  $D_2+dD_4$, with $d\ge 0$. This computation is very similar to the computation for $dD_4$ in the previous section. Thus, we omit details, but point out all the features that are different from the previous computation.

In order to compute the invariant $\langle D_i,pt \rangle_{0,2,D_2+dD_4}$ by localization, choose $[p_4]^\mathcal{T}=c_1^\mathcal{T} ( \calO^\mathcal{T} (D_2) ) \cdot c_1^\mathcal{T} ( \calO^\mathcal{T} (D_3) )$ as an equivariant lift for the insertion $[pt]$; note that  $D_2 \cap D_3 = p_4$.

Let $\beta = D_2 + dD_4$ and $F$ a fixed locus in $Q_{0,2}(\F_2,\beta)$. From the degree $\beta$, there is only one component of degree $D_2$ for a quasimap in $F$. The other components of degree $d'D_4$ with $d'\leq d$ must have their image on the $D_4$-curve in $\F_2$.
\begin{Lem}
\label{D3D2dD4}
    Let $\beta = D_2 + dD_4$ and $F$ a fixed locus in $Q_{0,2}(\F_2,\beta)$. Then, with the above choice of our equivariant lift for $[pt]$-insertion, for a quasimap in $F$, there is a component $C_0$ in the source curve $C$ satisfying
    \begin{enumerate}[label=(\roman*)]
        \item $C_0$ must be at one end in $C$ (a chain of $\PP^1$'s);
        \item the second marking must be on $C_0$;
        \item the image of $C_0$ lies on the $D_2$-curve;
        \item $C_0$ is of map-type.
    \end{enumerate}
    Hence, \[\langle D_3,pt \rangle_{0,2,D_2+dD_4}=0.\]
\end{Lem}
\begin{proof}
    Since we have a unique component of degree $D_2$, say $C_0$, (i) is clear.

    For (ii), due to our choice of the equivariant lift $[p_4]^\mathcal{T}$ for the $[pt]$-insertion, the second marking must be in $C_0$. Otherwise, there is no chance for the second marking to go to $p_4$.
    
    Note that $(D_2 \cdot D_\rho)_{\rho=1}^{4} = (0,0,1,1)$ implies that the image of $C_0$ lies on either $D_1$ or $D_2$. Since the second marking is in $C_0$, the image of $C_0$ must lie on the $D_2$-curve. Thus, (iii) is proved.

    Claim that $C_0$ cannot have any base points for (iv). If there is a base point, $C_0$ contracts to either $p_1$ or $p_4$. It must be $p_4$, since the second marking on $C_0$ must map to $p_4$. However, because the consecutive component of degree $d'D_4$ must land on the $D_4$-curve, contradiction occurs. Hence, $C_0$ must be of map-type.

    So far, the only component of degree $D_2$, which is of map-type and located at one end of the chain of $\PP^1$'s, has the second marking that goes to $p_4$. To verify $\langle D_3,pt \rangle_{0,2,D_2+dD_4}=0$, observe that the first marking is in a component of degree $d''D_4$ that is at the other end of the chain. This component maps to the $D_4$-curve. Recall that $D_3\cap D_4 = \varnothing$. Thus, there is no such a quasimap in $F$. Therefore, $\langle D_3,pt \rangle_{0,2,D_2+dD_4}=0$.
\end{proof}

From the relations in $H^*(\F_2)$, it suffices to compute
\[\langle D_1,pt \rangle_{0,2,D_2 + dD_4}.\]

In our chain graph expression for a fixed locus, we depict the component of degree $D_2$ by a vertical line. Fix this vertical line on the left. The chain graph of a general fixed locus is presented in Figure \ref{graphfixedlocus2}.


The $D_1$-insertion forces that the number of horizontal edges in the chain graph of a fixed locus must be odd so that the first marking maps to $D_1$.

To have a formula similar to \eqref{onestar} in Proposition \ref{formula}, we need to modify the notations and introduce more: for $k \in \{1,2\}$,
\begin{description}
    \item[$I_{k}^{m}$] the set of all interior vertices: 1) mapping to $p_k$, 2) not carrying any dashed half-edges, 3) not at the left end,
    \item[$I_{k}^{b}$] the set of all interior vertices: 1) mapping to $p_k$, 2) carrying a dashed half-edge, 3) not at the left end,
    \item[$I_2^{end}$] the set of all vertices at the right end: 1) mapping to $p_2$, 2) carrying a dashed half-edge,
    \item[$I_{D_2}^{m}$] the set of all interior vertices at the left end: 1) mapping to $p_1$, 2) not carrying any dashed half-edges,
    \item[$I_{D_2}^{b}$] the set of all interior vertices at the left end:  1) mapping to $p_1$, 2) carrying a dashed half-edge.
\end{description}
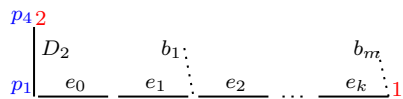
\begin{figure}[b]
    \centering
    \tikzset{every picture/.style={line width=0.75pt}} 
\begin{tikzpicture}[x=0.75pt,y=0.75pt,yscale=-1,xscale=1]
\draw    (8,300) -- (8,265) ;
\draw    (10,300) -- (45,300) ;
\draw    (50,300) -- (85,300) ;
\draw    (90,300) -- (125,300) ;
\draw    (150,300) -- (185,300) ;
\draw  [dash pattern={on 0.84pt off 2.51pt}]  (132,300) -- (144.6,300) ;
\draw  [dash pattern={on 0.84pt off 2.51pt}]  (83,275.6) -- (87.6,300) ;
\draw  [dash pattern={on 0.84pt off 2.51pt}]  (185,300) -- (179,274.25) ;
\draw (22,289.4) node [anchor=north west][inner sep=0.75pt]  [font=\tiny]  {$e_{0}$};
\draw (62,289.4) node [anchor=north west][inner sep=0.75pt]  [font=\tiny]  {$e_{1}$};
\draw (102,289.4) node [anchor=north west][inner sep=0.75pt]  [font=\tiny]  {$e_{2}$};
\draw (162,289.4) node [anchor=north west][inner sep=0.75pt]  [font=\tiny]  {$e_{k}$};
\draw (70,270.4) node [anchor=north west][inner sep=0.75pt]  [font=\tiny]  {$b_{1}$};
\draw (166,270.4) node [anchor=north west][inner sep=0.75pt]  [font=\tiny]  {$b_{m}$};
\draw (10,270.4) node [anchor=north west][inner sep=0.75pt]  [font=\tiny]  {$D_2$};
\draw (185,291.4) node [anchor=north west][inner sep=0.75pt]  [font=\tiny,color={rgb, 255:red, 255; green, 0; blue, 0 }  ,opacity=1 ]  {$1$};
\draw (-5,290.4) node [anchor=north west][inner sep=0.75pt]  [font=\tiny,color={rgb, 255:red, 0; green, 0; blue, 255 }  ,opacity=1 ]  {$p_1$};
\draw (-5,255.4) node [anchor=north west][inner sep=0.75pt]  [font=\tiny,color={rgb, 255:red, 0; green, 0; blue, 255 }  ,opacity=1 ]  {$p_4$};
\draw (7,255.4) node [anchor=north west][inner sep=0.75pt]  [font=\tiny,color={rgb, 255:red, 255; green, 0; blue, 0 }  ,opacity=1 ]  {$2$};
\end{tikzpicture}
    \caption{The chain graph of a general fixed locus for $\beta=D_2+dD_4$}
    \label{graphfixedlocus2}
\end{figure}
Collect the vertices regardless of where they go
\[I^m := I_{1}^{m} \sqcup I_{2}^{m}, \qquad I^b := I_{1}^{b} \sqcup I_{2}^{b}, \qquad I^B:=I^b \sqcup I_2^{end}.\]
Count them
\begin{align*}
&N_{k}^{m}:=|I_{k}^{m}|, \qquad N_{k}^{b}:=|I_{k}^{b}|, \qquad N_2^{end}:=|I_2^{end}|,~ \qquad~\\
    &N^{m}:=|I^{m}|,~ \qquad N^{b}:=|I^{b}|,\\
    &N_{D_2}^{m}:=|I_{D_2}^{m}|, \qquad N_{D_2}^{b}:=|I_{D_2}^{b}|.
\end{align*}
Note that $N_{D_2}^{m},~N_{D_2},~ N_2^{end} \in \{0,1\}$, and $N_{D_2}^{m} + N_{D_2}^{b}=1$.\\

One can have a formula for the inverse Euler class of the virtual normal bundle of a general fixed locus $F$ when $\beta = D_2 + dD_4$ through a similar way to have Proposition \ref{formula}. Let $e_{0}$ be the degree of the first horizontal edge from the left, and $b_{0}$ the degree of the dashed half-edge attached to the leftmost vertex, if it is in $I_{D_2}^{b}$.
\begin{Pro}
\label{formula2}
For a fixed locus $F$ in $Q_{0,2}(\F_2,D_2 + dD_4)$, is
\begin{align}
\label{twostar}
    &\frac{\big( 
    ev_{1}^{*}c_{1}^\mathcal{T}(\mathcal{O}^\mathcal{T}(D_1)) 
    ev_{2}^{*} \big( c_{1}^\mathcal{T}(\mathcal{O}^\mathcal{T}(D_2)) \cdot c_{1}^\mathcal{T}(\mathcal{O}^\mathcal{T}(D_3)) \big)
    \big)\big\lvert_{F}}{e^\mathcal{T}(N_{F}^{\text{vir}}) } \\ \nonumber
    &=  -\frac{W_1^2 V_1}{W_1V_1^2}
    \text{Cont}(VC) \text{Cont}_{m}(NS) \text{Cont}_{D_2}(NS) \\ \nonumber
    &\cdot \prod_{\text{edges}}\text{Cont}_{E}(e)\prod_{v \in I^{B}} \text{Cont}_{B}(v) \cdot \text{Cont}_{B}(D_2) , 
\end{align}
where 
\begin{align*}
    \text{Cont}(VC)&:= (-1)^{N_{2}^{m}+N_{2}^{end}+2N_{2}^{b}} V_1^{N_1^{m}+2N_1^{b}+N_{D_2}^{m} + 2N_{D_2}^{b} } V_2^{N_2^{m}+N_{2}^{end}+2N_2^{b}} W_1^{N^{m}+N_2^{end} + N_{D_2}^{m} + 2N_{D_2}^{b}+2N^{b}},\\
    \text{Cont}_{m}(NS)&:= \bigg( (-1)^{N_{2}^{m}} W_1^{N^{m}} \prod_{v\in I^{m}} \bigg (\frac{1}{e_v} + \frac{1}{e'_v} \bigg) \bigg)^{-1},\\
    \text{Cont}_{D_2}(NS)&:=\bigg(  (V_1+\frac{W_1}{e_{0}})^{N_{D_2}^{m}}   \bigg)^{-1},\\
    \text{Cont}_{E}(e)&:= \frac{e^{2e}\prod_{j=0}^{2e-2} \big( V_1 + \frac{1+j}{e}W_1 \big)}{(e!)^2W^{2e}(-1)^e}, \\
    \text{Cont}_{B}(v)&:= \frac{1}{b!} \int_{\overline{M}_{0,2|b}}
    \frac{ \frac{V_i}{W_i^2}\prod_{j=2}^{b}\frac{(V_i-2\Delta_{j})^2}{W_i+\Delta_j}}{ \bigg(\frac{W_i}{e} - \psi_{n}  \bigg)\bigg(\frac{W_i}{e'} - \psi_{n'}  \bigg)^{\epsilon(v)}},\\
    \text{Cont}_{B}(D_2)&:= \Bigg[ \frac{1}{b_{0}!} \int_{\overline{M}_{0,2|b_{0}}}
    \frac{ \frac{V_1}{W_1^2}\prod_{j=2}^{b_0}\frac{(V_1-2\Delta_{j})^2}{W_1+\Delta_j} }{ \bigg(\frac{W_1}{e_0} - \psi_{n}  \bigg)\bigg(V_1 - \psi_{n'}  \bigg) }  \Bigg]^{N_{D_2}^{b} }.\\
\end{align*}
\end{Pro}
\begin{proof}
    Observe that the weight of the numerator coming from the insertions $D_1$ and $pt$ is:
\[ \big( ev_{1}^{*}c_{1}^\mathcal{T}(\mathcal{O}^\mathcal{T}(D_1)) \big)  
 ev_{2}^{*} \big( c_{1}^\mathcal{T}(\mathcal{O}^\mathcal{T}(D_2)) \cdot c_{1}^\mathcal{T}(\mathcal{O}^\mathcal{T}(D_3)) \big)  \big\lvert_{F} = W^2V.\]

The edge contribution from the map-type component of degree $D_2$ is given by \[-\frac{1}{W_1V_1^2}\] as \eqref{term4}.

Also, we use the same formula when we compute \eqref{term1} to compute $\text{Cont}_{D_2}(NS)$.
\end{proof}

Besides the numerator, new factors appearing in formula \eqref{twostar} compared to formula \eqref{onestar} for the previous case are \[-\frac{1}{W_1V_1^2},~ N_{D_2}^{m},~N_{D_2}^{b},N_{D_2}^{end},~(V+\frac{W}{e_0})^{N_{D_2}^{m}},\]
which all are related to the component of degree $D_2$. We will pay attention to these in our analysis on the denominator of formula \eqref{twostar}.
\begin{Lem}
\label{denomWk2}
    The denominator of the formula \eqref{twostar} is of the form $W^N$ for some $N$.
\end{Lem}
\begin{proof}
The factor $\frac{W_1^2V_1}{W_1V_1^2}$ in formula \eqref{twostar} is $\frac{W_1}{V_1}$. We want to cancel the factor $\frac{1}{V} = \frac{1}{V_1}$. Note that we cannot just cancel this using the factor $V_1$ from either $\Cont_B(v)$ or $\Cont_B(D_2)$. The formal case is not possible because there is a map-type fixed locus that does not have this contribution term, or because there is a case where $v\in I^B$ might not map to $p_1$ so as to have $V_1=V$ factor. The latter case is not appropriate since the left corner vertex might not carry any dashed half-edges. Thus, we need to cancel the factor $\frac{1}{V}$ using a factor from some other places. Recall that $N_{D_2}^{m}+N_{D_2}^{b}=1$ for all fixed locus. Therefore, we can cancel $\frac{1}{V}$ by $V_1^{N_{D_2}^{m}+N_{D_2}^{b}}$ in $\Cont(VC)$.

The next factor we need to consider is $(V_1+\frac{W_1}{e_{0}})^{N_{D_2}^{m}}$ in $\Cont_{D_2}(NS)$. Observe that the horizontal edge at the leftmost end is the map-type component of degree $e_0=e$. Thus, with $j=0$, there is a factor $(V_1+\frac{1+j}{e}W_1)$ in $\Cont_E(e)$ to cancel the factor $(V_1+\frac{W_1}{e_{0}})^{N_{D_2}^{m}}$ in $\Cont_{D_2}(NS)$.

Last, the factor $\frac{1}{V_1}$ from $\frac{1}{V_1 - \psi_{n'}} = \frac{1}{V_1}\sum_{l}(\psi_{n'}/V)^{l}$ is cancelled by the factor $V_1$ coming from the numerator in $\Cont_{B}(D_2)$.

Hence, the denominator of formula \eqref{twostar} is $W^N$ for some $N$.
\end{proof}
Since we figured out the form of the denominator in formula \eqref{twostar}, we can decide necessary fixed loci whose contribution is nonzero. 
\begin{Cor}
\label{correducedloci2}
The general graphs of necessary fixed loci for $\beta=D_2+dD_4$ are in Figure \ref{reducedloci2}, where $b+e=d$. We denote them as $F'_b$.
\end{Cor}
\begin{figure}[b]
\tikzset{every picture/.style={line width=0.75pt}} 
\begin{tikzpicture}[x=0.75pt,y=0.75pt,yscale=-1,xscale=1]
\draw    (60,109) -- (97,109) ;
\draw    (59,109) -- (59,79) ;
\draw  [dash pattern={on 0.84pt off 2.51pt}]  (85,90) -- (100,109) ;
\draw (75,100) node [anchor=north west][inner sep=0.75pt]  [font=\tiny]  {$e$};
\draw (88,84) node [anchor=north west][inner sep=0.75pt]  [font=\tiny]  {$b$};
\end{tikzpicture}
\caption{Necessary fixed loci for $\beta=D_2+dD_4$}
\label{reducedloci2}
\end{figure}
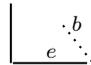
\begin{proof}
Consider the factor
\[\frac{W^2V}{WV^2}V^{N_1^{m}+2N_1^{b}+N_{D_2}^{m} + 2N_{D_2}^{b}}\]
in formula \eqref{twostar} has positive exponent, when $N_{D_2}^{b}=1$. Because of Lemma \ref{denomWk2}, the fixed loci with $N_{D_2}^{b}=1$ has zero contribution to the final answer. Thus, we will only consider the case $N_{D_2}^{m}=1$ from $N_{D_2}^{m}+N_{D_2}^{b}=1$.
Meanwhile, the exponents of $V_1$, $N_1^{m}$, in $\Cont(VC)$ must be zero to give nonzero contribution. Hence, the only possible fixed loci with nonzero contribution have \[
N_{D_2}^{b}=0,~N_{D_2}^{m}=1,~N_1^{m}=0,~2N_1^{b}=0,~N_{2}^{end}=1.\]
\end{proof}

For $\beta = D_2 + dD_4$ with $b+e= d$, let $F$ be a necessary fixed locus depicted in Figure \ref{reducedloci2}. Applying the contents in Sections \ref{comb1}, \ref{comb2}, and \ref{comb3} with setting $V=0$, one can have the following simplification of formula \eqref{twostar}:
\begin{align}
\label{twostarsimp}
    (-1)^{e-1}\frac{e(2e-1)!}{(e!)^2}\bigg(\frac{4^b e}{b}\binom{b+e-1}{b-1} \bigg)^{N_{2}^{end}} =  \frac{(-1)^{d-b-1}4^b}{2}\binom{2(d-b)}{d-b}\binom{d-1}{b}.
\end{align}
We give a computation of the invariants $\langle D_1,pt \rangle_{0,2,D_2 + dD_4}$. 
\begin{Pro}
\label{computeD2dD4}
All 2-pointed degree $D_2 + dD_4$ quasimap invariants of $\F_2$ are given by
\begin{align*}
    \langle D_1,pt \rangle_{0,2,D_2 + dD_4}&=\langle D_2,pt \rangle_{0,2,D_2 + dD_4} = \frac{1}{2(2d-1)}\binom{2d}{d},
    \\
    \langle D_3,pt \rangle_{0,2,D_2 + dD_4}&=0,\quad \langle D_4,pt \rangle_{0,2,D_2 + dD_4} = -\frac{1}{(2d-1)}\binom{2d}{d}.
\end{align*}
\end{Pro}
\begin{proof}
Apply the Atiyah--Bott localization theorem to $\langle D_1,pt \rangle_{0,2,D_2 + dD_4}$, and use \eqref{twostarsimp}. Then, similar argument in the proof of Proposition \ref{computedD4} gives
    \[ \sum_{b=0}^{d-1}\frac{(-1)^{d-b-1}4^b}{2}\binom{2(d-b)}{d-b}\binom{d-1}{b} = -\frac{4^d}{2}\sum_{b=0}^{d-1}\binom{d-1}{d}\binom{-1/2}{d-b}.\]
Thus, it is enough to show \[\sum_{b=0}^{d-1}\binom{d-1}{b}\binom{-1/2}{d-b} = -\frac{4^{-d}}{(2d-1)}\binom{2d}{d}.\] The coefficient of $x^d$ in $(1+x)^{d-1-1/2}$ is $\sum_{b=0}^{d-1}\binom{d-1}{b}\binom{-1/2}{d-b}$. We use induction. Let $a_{d+1}$ be the coefficient of $x^{d+1}$ in $(1+x)^{d-1/2}$. Then, applying $\frac{d}{dx}$ and the induction hypothesis give us
\[a_{d+1} = -\frac{4^{-d}}{2d-1}\frac{d-1/2}{d+1}\binom{2d}{d} =-\frac{4^{-(d+1)}(2d+2)(2d+1)}{(2d+1)(d+1)^2}\binom{2d}{d}=-\frac{4^{-(d+1)}}{(2d+1)}\binom{2(d+1)}{d+1}. \]

Lemma \ref{D3D2dD4} computes $\langle D_3,pt \rangle_{0,2,D_2 + dD_4}$. The relation allows us to compute $\langle D_1,pt \rangle_{0,2,D_2 + dD_4}$ and $\langle D_4,pt \rangle_{0,2,D_2 + dD_4}$.
\end{proof}

From the relation $D_3 = 2D_1+D_4$, one can derive \[\langle D_4,pt \rangle_{0,2,D_2 + dD_4}=-\frac{1}{2d-1}\binom{2d}{d}.\]

Propositions \ref{computedD4} and \ref{computeD2dD4} are necessary ingredients to prove the main Theorem \ref{Thm} in this paper.
\begin{proof}[proof of Theorem \ref{Thm}]
    We show the computation of $\sigma_2 \star 1$, and the rest can be computed in a very similar way. Let $q_2:= q^{D_2}$ and $q_4:= q^{D_4}$. Recall that $f(z):=\sum_{d\geq 1}\binom{2d}{d}z^d = \frac{1}{\sqrt{1-4z}}-1$. Applying Proposition \ref{computedD4}, the divisor equation \eqref{diveqn}, and the Point mapping axiom and the Degree axiom in \cite[\S7.3]{coxmirror} give us
    \begin{align*}
        \sigma_2 \star 1 &:= \sum_i \sum_{\beta \in \text{Eff}} q^\beta \langle 1, T_i \mid \sigma_2 \rangle_{0,2|1,\beta} T^i\\
        &=\sum_{d \geq 0} q_4^{d}\langle 1, D_2 \mid \sigma_2 \rangle_{0,2|1,dD_4} (2D_2 + D_4) + \sum_{d \geq 0} q_4^{d}\langle 1, D_4 \mid \sigma_2 \rangle_{0,2|1,dD_4} D_2\\
        &=(D_2\cdot D_2)(2D_2 + D_4)+\sum_{d \geq 1} q_4^{d}(D_2 \cdot dD_4)\frac{-1}{2d}\binom{2d}{d} (2D_2 + D_4)\\
        &\quad +(D_2\cdot D_4)D_2+\sum_{d \geq 1} q_4^{d}(D_2 \cdot dD_4)\frac{1}{d}\binom{2d}{d} D_2\\
        &=D_2 -\frac 12 f(q_4)D_4.
    \end{align*}
    
\end{proof}

\section{The Batyrev Quantum Ring for $\F_2$}
\label{sec_conj}
For a smooth projective toric variety $X_\Sigma$, Batyrev defined in \cite{batyrev} a ring  from the data of the fan $\Sigma$.
In \cite[\S1.1.8]{lucathesis}, the author pointed out that the quantum deformation given by 3-pointed quasimap invariants is not the same as the Batyrev ring of $\F_2$ due to the failure of the divisor equation. In this section, we show that the quantum module structure of $\F_2$ in Theorem \ref{module} using $2|1$-pointed quasimap invariants agrees with the Batyrev ring of $\F_2$ realized as a natural module over the ring $\C[[q_2,q_4]]$.

\subsection{The Batyrev ring of $X_\Sigma$}\label{subsecbat} Let $v_1,\ldots, v_s \in N \cap \Sigma(1)$ be primitive integral generators for the rays. There are two ideals in $\C[x_1,\ldots,x_s]$. The first ideal is given by
\begin{align*}
    P_\Sigma := \big\langle \sum_{i=1}^{s}\langle m, v_i \rangle x_i \mid m \in M \big\rangle.
\end{align*}
For a primitive collection $P=\{v_{i_1},\ldots, v_{i_k}\}$, we have the relation
\begin{equation}
\label{relrays}
    v_{i_1}+\cdots+v_{i_k} = c_{1}v_{j_1}+\cdots+c_{l}v_{j_l},
\end{equation}
where $v_{j_1},\ldots,v_{j_l}$ are the generators of $\sigma \in \Sigma$ such that $v_{i_1}+\cdots+v_{i_k} \in \sigma$ and $c_1,\ldots, c_l \geq 0$. Using the dual of the exact sequence \eqref{seqclassgroup}, the relation \eqref{relrays} gives rise to a class $\beta_P$ in $H_2(X,\Z)$. This class $\beta_P$ is effective \cite[Thm 2.15]{batyreveff}, \cite[Example 8.1.2.2]{coxmirror}. Then, the second ideal is defined by
\begin{align*}
    SP_{\Sigma} := \langle x_{i_1}\cdots x_{i_k}-q^{\beta_P} x_{j_1}^{c_1}\cdots x_{j_l}^{c_l}\mid P:\text{ primitive collection} \rangle 
\end{align*}
This ideal is called the \textit{quantum Stanley Reisner ideal}.

Using these ideals, Batyrev defined the following ring
\[\text{Bat}H^*(X):= \C[x_1,\ldots, x_s] / (P_\Sigma + SP_{\Sigma}).\]

When $X_\Sigma$ is Fano, the quantum cohomology ring of $X_\Sigma$ coincides with the Batyrev ring; however, this is false when $X_\Sigma$ is not Fano, but \textit{semipositive}, i.e., the anticanonical divisor is nef. The Hirzebruch surface of type 2, $\mathbb{F}_2$, exactly shows the failure, see \cite[Example 11.2.5.2]{coxmirror}.

Since the effective cone of $\F_2$ is generated by $D_2$ and $D_4$, one can write down all $\beta_P$ as the nonnegative linear combination of $D_2$ and $D_4$. The Batyrev ring for $\F_2$ can be written as 
\[\text{Bat}H^*(\F_2)=\mathbb{C}[x_2,x_4]/\langle x_2^2-q_4x_4^2, 
(2x_2+x_4)x_4 -q_2  \rangle.\]

The quantum $H^*((\C^*)^2)$-module of $\F_2$ has the following relations.
\begin{Lem}
\label{relinquantum}
The following relations hold in the quantum $H^*((\C^*)^2)$-module of $\F_2$:
\begin{equation*}
    (2\sigma_2 + \sigma_4) \star (\sigma_4 \star 1) =q_2, \quad \sigma_2 \star (\sigma_2  \star 1) = q_4 \sigma_4 \star (\sigma_4  \star 1).
\end{equation*}
\end{Lem}
\begin{proof}
Let $f(z)=\sum_{d\geq 1}\binom{2d}{d}z^d = \frac{1}{\sqrt{1-4z}}-1$. Then, one can verify
\begin{align}
\label{relfandz}
(1+f)^2(1-4z)=1, \\
4z(1+f)^2=f(2+f).  \label{relfandz2}  
\end{align}
Observe the following:
\begin{align*}
    (2\sigma_2 + \sigma_4) \star (\sigma_4 \star 1) &= (1+f)(2\sigma_2 + \sigma_4) \star D_4 \\
    &=q_2(1+f)(-4q_4(1+f) +(1+f) )\\
    &=q_2(1+f)^2(-4q_4 +1 ) \stackrel{\eqref{relfandz}}{=} q_2,
\end{align*}
and 
\begin{align*}
    \sigma_2 \star (\sigma_2  \star 1) &= \sigma_2 \star (D_2-\frac 12 fD_4)\\
    &= \bigg( q_2q_4(1+f) -\frac{1}{2}fpt \bigg) -\frac 12 f \bigg(-2q_2q_4(1+f) + (1+f)pt\bigg)\\
    &= q_2q_4 (1+f)^2 -\frac 12 pt f(2+f) \stackrel{\eqref{relfandz2}}{=} q_2q_4 (1+f)^2 -pt 2  q_4(1+f)^2\\
    &= q_4(1+f)^2 (q_2 - 2 pt),\\
    \sigma_4 \star (\sigma_4  \star 1) &= (1+f)\sigma_4 \star D_4\\
    &= (1+f) \bigg( q_2(1 + f)-2(1+f)pt \bigg) \\
    &=  (1+f)^2 ( q_2-2pt).
\end{align*}
Here, $f$ stands for $f(q_4)$. Thus, we have \[ \sigma_2 \star (\sigma_2  \star 1) = q_4 \sigma_4 \star (\sigma_4  \star 1).\]
\end{proof}
One can see that the two relations in Lemma \ref{relinquantum} are matching up with the relations in the Batyrev ring $\text{Bat}H^*(\F_2)$ by viewing $\sigma_2 \star 1$ and $\sigma_4 \star 1$ as $x_2$ and $x_4$, respectively.
\begin{Lem}
\label{nakayama}
    The Batyrev ring $\textup{Bat}H^*(\F_2)$ is generated by $1,x_2,x_4$ and $x_2x_4$ over $\C[[q_2,q_4]]$.
\end{Lem}
\begin{proof}
    Applying Nakayama's lemma to the local ring $\C[[q_2,q_4]]$ with the unique maximal ideal $( q_2,q_4 )$, one can show that the set $\{1,x_2,x_4,x_2x_4\}$ generates $\textup{Bat}H^*(\F_2)$ over $\C[[q_2,q_4]]$.
\end{proof}
\begin{Rem}
    The above proof of Lemma \ref{nakayama} does not give us a concrete way to write down an element in $\textup{Bat}H^*(\F_2)$ with respect to the generating set $\{1,x_2,x_4,x_2x_4\}$. It is worthwhile to observe, for instance, how $x_2^2x_4$ can be written in terms of $1,x_2,x_4$, and $x_2x_4$, since it requires the equation \eqref{relfandz}, which comes from the series $f(z)=\sum_{d\geq 1} \binom{2d}{d} z^d $ that contains some 2-pointed quasimap invariants of $\F_2$.
    \begin{align*}
        x_2^2 x_4 &= q_4 x_4^3 = x_4q_4(q_2 - 2x_2 x_4 ) = q_2q_4x_4 - 2q_4x_2 x_4^2 \\
        &=q_2q_4x_4 - 2q_4x_2 (q_2-2x_2x_4) = q_2q_4x_4 - 2q_2q_4x_2 +4 q_4 x_2^2x_4\\
        (1-4q_4)x_2^2 x_4 &= q_2q_4(x_4-2x_2)\\
        x_2^2 x_4 &= q_2q_4(1+f(q_4))^2(x_4-2x_2).
    \end{align*}
    This shows that finding generating sets over $\C[[q_2,q_4]]$ using the relations in the quantum cohomology ring involves with some generating series whose coefficients are given by 2-pointed quasimap invariants.
\end{Rem}

A natural $H^*((\mathbb{C}^*)^2 ) = \C[\sigma_2,\sigma_4]$-module structure on $\textup{Bat}H^*(\F_2)$ is given by
\[\sigma_2 \cdot x_2^a x_4^b = x_2^{a+1} x_4^b, \quad \sigma_4 \cdot x_2^a x_4^b = x_2^{a} x_4^{b+1}.\]
For convenience, we call this the \textit{Batyrev module}. 
\begin{Pro}
    The $\C[\sigma_2,\sigma_4]$-module $\textup{Bat}H^*(\F_2)$ is isomorphic to the quantum $H^*( (\mathbb{C}^*)^2 )$-module of $\F_2$ from Theorem \ref{Thm}.
\end{Pro}
\begin{proof}
Denote the generating sets by $\alpha:=\{1,x_2,x_4,x_2x_4\}$ and $\mu:= \{1,D_2,D_4,pt\}$ for the Batyrev module and the quantum module, respectively. Using Lemma \ref{nakayama}, one can represent the action of $\sigma_i$ as matrices using $\alpha$ and $\mu$, respectively, say $[\sigma_i]_{\alpha}^{\alpha}$ and $[\sigma_i]_{\mu}^{\mu}$. 

Define a linear function $\phi$ in the following way:
\[1 \mapsto 1, \quad x_2 \mapsto \sigma_2 \star 1, \quad x_4 \mapsto \sigma_4 \star 1, \quad x_2x_4 \mapsto (\sigma_2\sigma_4) \star 1.\]
This map extends linearly over $\C[[q_2,q_4]]$. Then, the matrix presentation of $\phi$ over $\C[[q_2,q_4]]$ is the following:
\begin{equation*}
[\phi]_{\mu}^{\alpha}:= \begin{pmatrix}
    1&0&0&-2(1+f)^2q_2q_4 \\
    0&1&0&0 \\
    0&(-1/2)f&1+f&0 \\
    0&0&0&(1+f)^2 
\end{pmatrix}, 
\end{equation*}
where $f$ is denoted for $f(q_4)$. One can check
\[ [\phi]_{\mu}^{\alpha} [\sigma_i]_{\alpha}^{\alpha}=[\sigma_i]_{\mu}^{\mu}[\phi]_{\mu}^{\alpha},\]
which shows that $\phi$ is $\C[\sigma_2,\sigma_4]$-linear over $\C[[q_2,q_4]]$. From Lemma \ref{relinquantum}, there exists a well-defined induced map from the Batyrev module to the quantum module. The determinant of $[\phi]_{\mu}^{\alpha}$ is $(1+f)^3$. Since by equation \eqref{relfandz} with $z=q_4$, $(1+f)^3$
is an invertible element in $\C[[q_2,q_4]]$, so that $\phi$ is an $\C[\sigma_2,\sigma_4]$-isomorphism over $\C[[q_2,q_4]]$.
\end{proof}

{
\bibliographystyle{plain}
\bibliography{References.bib}
}
\end{document}